\documentclass[10pt,fullpage]{article}
\usepackage{enumerate}
\usepackage{pgf}
\usepackage{tikz}
\usetikzlibrary{arrows}
\usepackage[utf8]{inputenc}

\newtheorem{theorem}{Theorem}
\newtheorem{lemma}{Lemma}

\newtheorem{proposition}{Proposition}

\begin{document}
\title{On mild contours in ray categories}

\author{Klaus Bongartz\\Universität Wuppertal\\Germany\thanks{E-mail:bongartz@math.uni-wuppertal.de}}
\maketitle
\begin{abstract}
 We generalize and refine the structure and disjointness theorems for non-deep contours obtained in the fundamental article 'Multiplicative bases and representation-finite algebras'. In particular we show that these contours do not occur in minimal representation-infinite algebras.
\end{abstract}

\section*{Introduction}

The article 'Multiplicative bases and representation-finite algebras'  written by Bautista, Gabriel, Roiter and Salmer\'{o}n contains many beautiful results and ideas and it had some nice consequences: the proof by Bautista et alii of the second Brauer Thrall conjecture, Fischbachers result on the universal cover of a ray category implying that my criterion for finite representation type always works and recently  the assertion that there are no gaps in the lengths of the indecomposables.

 The article on multiplicative bases can be divided into two main parts: a local one dealing with small parts of the algebra and a global one dealing with topological and cohomological questions. The local part contains the surprisingly simple structure and disjointness theorems for non-deep contours in mild categories i.e. distributive categories such that any proper quotient is locally representation-finite. These theorems are crucial to show the existence of a semi-multiplicative basis which is then transformed into a multiplicative basis by using a vanishing theorem on cohomology proven in the global part.

 In the present paper we generalize and refine the main results of the local part. First, by taking a more direct approach to non-deep contours we obtain a very natural proof of the structure theorem. It is independent of any lists and  also much shorter especially in the most complicated case of diamonds. Second, we show that two different non-deep contours do not even share a point, and finally we prove that non-deep contours only occur in the locally representation-finite case. An a priori  proof of this last fact would give an almost conceptual proof of the second Brauer-Thrall conjecture.

\section{Reminder on ray categories and statement of the results}
To state our results precisely and to fix the notations we recall some basic facts about ray categories which can all be found in  chapter 13 of \cite{Buch}.

 A finite dimensional module  $M$ over a locally bounded $k$-category $A$ is just a covariant $k$-linear functor from $A$ to the category of $k$-vectorspaces such that the sum of the dimensions of all $M(x), x \epsilon A$, is finite. $A$ is distributive  if all endomorphism algebras $A(x,x)$ are uniserial and all homomorphism spaces $A(x,y)$ are cyclic as an $A(x,x)$ right module or as an $A(y,y)$ left module. The product $A(x,x)^{\ast}\times A(y,y)^{\ast}$ of the two automorphism groups acts on $A(x,y)$ and the orbit of a morphism is the corresponding ray. These rays are the morphisms of the ray category $\vec{A}$ attached to $A$ which is a kind of combinatorial skeleton underlying $A$. The properties of $\vec{A}$ are subsumed in the following axioms that define the abstract notion of a ray category $P$:
\begin{enumerate}
 \item The objects form a set and they are pairwise not isomorphic.
\item There is a family of zero-morphisms $0_{xy}:x \rightarrow y $, $x,y \epsilon P$, satisfying $\mu 0=0=0\nu$ whenever the composition is defined.
\item For each $x \epsilon P$, $P(x,y)=\{0\}$ and $P(y,x)=\{0\}$ for almost all $y \epsilon P$.
\item For each $x$ one has $P(x,x)=\{ id_{x},\sigma,\ldots ,\sigma^{n-1}\neq 0=\sigma^{n}\}$. Here $n$ depends on $x$.
\item For each $x,y$, the set $P(x,y)$ is cyclic under the action of $P(x,x)$ or $P(y,y)$.
\item If $\kappa, \lambda,\mu,\nu$ are morphisms with $\lambda \mu\kappa= \lambda \nu \kappa\neq 0$ then $\mu= \nu$.
 \end{enumerate}
We denote the generator $\sigma$ in d) by $\pi(x)$ and the generator in e) by $\pi(x,y)$. 
It follows from property e) that the 'subbimodule' of $P(x,y)$ generated by a morphism $\mu$ is cyclic over $P(y,y)$ or over $P(x,x)$ or over both. The morphism is then called a transit or cotransit or bitransit morphism. If $\mu\neq 0$ is annihilated by $\pi(x)$ and by $\pi(y)$ it is called deep. The generator $\pi(x,y)$ has the same transit-properties as any non-deep morphism from $x$ to $y$.

 Starting with an abstract ray category $P$ one constructs in a natural way its linearization $k(P)$, which is a locally bounded distributive $k$-category having the original category $P$ as the associated ray category $\vec{k(P)}$.  We say that $P$ is  ( locally ) representation-finite or minimal representation-infinite if $k(P)$ is so, and this is independent of the field by the well-known criterion on finite representation type
\cite[theorem 14.7]{Buch}.
To study a ray category $P$ and its universal cover $\tilde{P}$ we look at the quiver $Q_{P}$ of $P$. Its points are the objects of $P$ and its arrows the irreducible morphisms in $P$, i.e. those non-zero morphisms that cannot be written as a product of two morphisms different from identities. Therefore we do not  distinguish  between 'arrows' and 'irreducible morphisms' when working in a ray category. The path category $\mathcal{P}Q_{P}$ has the points of $Q_{P}$ as objects and the paths in $Q_{P}$ as non-zero morphisms, to which we add formal zero-morphisms. There is a canonical full functor $\vec{}:\mathcal{P}Q_{P} \longrightarrow P$ from the path category to $P$ which is the 'identity' on objects, arrows and zero-morphisms. 

Two paths in $Q_{P}$ are interlaced if they belong to the transitive closure of the relation $R$ given by $(v,w)\epsilon R$ iff $v=pv'q, w=pw'q$ and $\vec{v}'=\vec{w}'\neq 0$ where $p$ and $q$ are not both identities. A contour $C$  of $P$ from $x$ to $y$ is a set $C= \{v,w\}$ of non-interlaced paths from $x$ to $y$ with $\vec{v}=\vec{w} \neq 0$. The contour is non-deep if $\vec{v}$ is so. Throughout this article we denote by $P(C)$ the full subcategory of $P$ supported by the points in $v$ and $w$ and by $Q(C)$ the quiver of $P(C)$ which is in general not a subquiver of $Q_{P}$.

Figure 1.1 describes three  ( families of ) ray categories by quiver and relations. Each of these contains a non-deep contour from $x$ to $y$ and $\pi(x,y)$ is always bitransit. For obvious reasons the contours $C$ as well as the  categories $P(C)$ are called  penny-farthings, dumb-bells and diamonds respectively.

\setlength{\unitlength}{1cm}
\begin{picture}(15,8)

\put(8,4){\circle*{0.1}}
\put(7,2){\circle*{0.1}}
\put(5,1){\circle*{0.1}}
\put(0,4){\circle*{0.1}}
\put(3,1){\circle*{0.1}}
\put(1,2){\circle*{0.1}}
\put(1,6){\circle*{0.1}}
\put(3,7){\circle*{0.1}}
\put(5,7){\circle*{0.1}}
\put(7,6){\circle*{0.1}}
\put(8.5,4){\circle{1}}
\put(7,6){\vector(1,-2){1}}
\put(8,4){\vector(-1,-2){1}}
\put(8.0,3.9){\vector(0,1){0.1}}
\put(9.2,4){$\rho$}
\put(5,7){\vector(2,-1){2}}
\put(7,2){\vector(-2,-1){2}}
\put(5,1){\vector(-2,0){2}}
\put(3,7){\vector(2,0){2}}
\put(3,1){\vector(-2,1){2}}
\put(1,6){\vector(2,1){2}}
\multiput(0,4)(0.1,0.2){10}{\circle*{0.01}}
\multiput(0,4)(0.1,-0.2){10}{\circle*{0.01}}

\put(7,5){$\alpha_{n}$}
\put(7,3){$\alpha_{1}$}
\put(5.8,3.9){$x=y=x_{0}$}
\put(6.5,2){$x_{1}$}

\put(5.8,6){$x_{n-1}$}
\put(8,3){$n\geq 1,\alpha_{1}\alpha_{n}=0$}
\put(8,2.5){v=$\alpha_{n} \ldots \alpha_{1}, w =\rho^{2}$}
\put(7.5,2){$ 0= \alpha_{e(i)}\ldots \alpha_{1}\rho\alpha_{n}\ldots \alpha_{i+1}$}
\put(7,1.5){$e:\{1, \ldots ,n-1\} \rightarrow \{1,\ldots ,n\}$}
\put(7,1){$e$ non-decreasing}

\end{picture}

\setlength{\unitlength}{0.5cm}
\begin{picture}(20,7)

\put(11,4){\circle*{0.2}}
\put(15,7){\circle*{0.2}}
\put(15,1){\circle*{0.2}}
\put(19,4){\circle*{0.2}}
\put(11,4.2){\vector(4,3){3.8}}
\put(15,6.8){\vector(-4,-3){3.8}}
\put(11,4){\vector(4,-3){3.8}}
\put(15,1){\vector(4,3){3.8}}
\put(18.8,4){\vector(-4,3){3.8}}
\put(15.2,7){\vector(4,-3){3.8}}

\put(10,4){$x$}
\put(19.5,4){$y$}
\put(15.5,7.2){$z$}
\put(15.5,0.7){$t$}

\put(12,5.5){$\gamma$}
\put(14,5.5){$\lambda$}
\put(16,5.5){$\kappa$}
\put(17.5,5.5){$\alpha$}

\put(12,2.5){$\delta$}
\put(17.5,2.5){$\beta$}

\put(11,0){$v=\beta\delta, w = \alpha \gamma$,}
\put(16,0){$\lambda \kappa =0,$}
\put(19,0){$\kappa \alpha= \gamma \lambda$}

\put(1,4){\circle*{0.1}}
\put(3,4){\circle*{0.1}}
\put(0.5,4){\circle{1}}
\put(3.5,4){\circle{1}}
\put(1,4){\vector(2,0){2}}
\put(1,4.1){\vector(0,-1){0.1}}
\put(3,4.1){\vector(0,-1){0.1}}
\put(0.5,3.9){$x$}
\put(3.3,3.9){$y$}
\put(1,4.5){$\lambda$}
\put(2,4.5){$\mu$}
\put(3,4.5){$\rho$}
\put(1,3){$v=\mu\lambda,w = \rho \mu$}
\put(1,2){$\lambda^{r} =0=\rho^{s}$}
\put(0,1){ $min \{r,s\} =3, max \{r,s\} \leq$ 5}
\put(9,-2){figure 1.1}
\end{picture}\vspace{2cm}

The universal covers are easy to determine in all three cases. The fundamental group is always ${\bf
Z }$ and the three categories $P(C)$ are representation-finite by the criterion or by a direct calculation of the Auslander-Reiten-quivers using the method from \cite[section 14.4]{Buch}.

 Let $C=\{v,w\}$ be a non-deep contour in $P$ from $x$ to $y$ and suppose that $\pi(x,y)$ is transit. A full subcategory $D$ of $P/\pi(y)\vec{v}$  is called decisive   for $C$ if $D$ contains $x$ and $y$ and is contained in the union of the supports of the projective $P_{x}$ to $x$ and the injective $I_{y}$ to $y$. Then $C$ is called mild if $D$ is representation-finite for all full subcategories $D$ of $P/\pi(y)\vec{v}$ that have at most 4 points and are decisive for $C$. In case $\pi(x,y)$ is cotransit the mildness of $C$ is defined dually. 

Our first result reads as follows:
\begin{theorem}\label{structure} Any mild contour  is a dumb-bell, a penny-farthing or a diamond.
 
\end{theorem}

This generalizes considerably the structure theorem in  \cite{BGRS} and also the more precise statement in part a) of theorem 13.12 in \cite{Buch}. Namely, we only consider subcategories with at most 4 points in only one quotient.

To formulate our disjointness-result let $C$ be a mild contour from $x$ to $y$ as above and let $C'$ be a contour from $x'$ to $y'$. A full subcategory $D$ of $P/\pi(y)\vec{v}$ is decisive for the pair $(C,C')$ if it contains $x$ and $y$ and is contained in the union of the supports of $P_{x}$ and $P_{x'}$. Note that this definition is not symmetric.

\begin{theorem}\label{disjoint} 
 Let $C$ and $C'$ be two mild contours.
\begin{enumerate}
 \item Suppose $D$ is representation-finite for all full subcategories $D$ of $P/\pi(y)\vec{v}$ that are decisive for $(C,C')$ and have at most 5 points.
Then no arrow of $v$ or $w$ occurs in $v'$ or $w'$.
\item Suppose $D$ is representation-finite for all full subcategories $D$ of $P/\pi(y)\vec{v}$ that are decisive for $(C,C')$ and have at most 6 points.
Then no point  of $Q(C)$ occurs in $Q(C')$.
\end{enumerate}

\end{theorem}

 This refines the disjointness theorem in \cite{BGRS} and also its more precise version in section 13.12 of \cite{Buch}.

 The next result was known before only for penny-farthings by \cite{Standard}.

\begin{theorem}\label{minimal}
 If $P$ is minimal representation-infinite it contains only deep contours.
\end{theorem}

We give two independent proofs for this result. The first one uses the structure theorem obtained before and  it consists in a detailed study of the way a mild contour is embedded into the whole ray category. The second proof does not depend on the structure theorem but it uses the fact shown in  \cite{gaps} that a minimal representation-infinite ray category has an interval-finite universal cover and a free fundamental group. Thus one can apply the finiteness criterion and one finds directly with the help of the Bongartz-Happel-Vossieck list that only dumb-bells and penny-farthings with $n=2$ have to be excluded. In the following we abbreviate this list by BHV-list and we refer to the numbering from \cite[10.7]{Buch}.

Our main working tools are cleaving diagrams as in \cite{BGRS}. Therefore we recall the relevant facts and notations.
A diagram $D$ in $P$ is just a covariant functor $F:D \rightarrow P$ from another ray category $D$ to $P$. Then $F$ ( or also $D$ ) is called cleaving  iff it satisfies the following two conditions and their duals: a) $F\mu =0$ iff $\mu=0$;
b) If $\alpha \epsilon D(x,y)$ is irreducible and $F\mu:Fx \rightarrow Fz$ factors through $F\alpha$ then $\mu$ factors already through $\alpha$.
In practice, the following equivalent conditions from \cite{gaps} are easier to verify: a) $F\mu=0$ iff $\mu=0$; b) No irreducible morphism is mapped to an identity; c) For any two irreducible morphisms $\alpha:x \rightarrow y$  and $\beta:x \rightarrow z$ in $D$ and each non-zero morphism  $\mu=\mu '\beta:x \rightarrow t$ with $\gamma\mu=0$ for all arrows $\gamma$ the image $F\mu$ does not factor through $F\alpha$ provided $\mu$ does not factor through $\alpha$; d) The dual of c).

The key fact about cleaving functors is that $P$ is not ( locally ) representation-finite if $D$ is not. In this article $D$ will always be given by its quiver $Q_{D}$, that has no oriented cycles, and some relations. Two  paths between the same points give always the same morphism, and zero relations are written down explicitely. The cleaving functor is then defined by drawing the quiver of $D$ with relations and by writing the morphism $F\alpha$ in $P$ close to each arrow $\alpha$. 

For  example let $D$ be the ray category with the natural numbers as objects and with arrows $2n \leftarrow 2n+1$ and $2n+1 \rightarrow 2n+2$ for all $n$. Then a cleaving functor from $D$ to $P$ is called a zigzag  and $P$ is said to contain a zigzag. A functor from $D$ to $P$ is just an infinite sequence of morphisms $(\sigma_{1},\rho_{1},\sigma_{2},\rho_{2},\ldots   )$ in $P$ such that $\rho_{i}$ and $\sigma_{i}$ always have common domain and $\rho_{i}$ and $\sigma_{i+1}$ common codomain.
 The  functor is cleaving iff none of the equations $\sigma_{i}=\xi \rho_{i}$,$\xi \sigma_{i}= \rho_{i}$,$\sigma_{i+1}\xi =\rho_{i}$ or $\sigma_{i+1}= \rho_{i} \xi$ has a solution. A crown in $P$  is a zig-zag that becomes periodic after n steps, i.e. one has  $\sigma_{i}=\sigma_{n+i}$ and $\rho_{i}=\rho_{n+i}$.   We denote such a crown by $\tilde{A}(\sigma_{1},\rho_{1},\sigma_{2}, \ldots ,\rho_{n})$. \vspace{0.5cm}

In fact, the three theorems  can be made more precise in the spirit of \cite[13.10,13,11,13.15]{Buch}. Namely, we do not really use that the considered categories are representation-finite but only the fact that apart from extended Dynkin quivers certain subcategories belonging to some frames of the BHV-list do not occur as cleaving diagrams. For instance, in the proof of the structure theorem we have to exclude the frames with numbers 11,14 and 21 to which one has to add number 20 for the proof of the disjointness theorem. From this point of view the most complicated part is the analysis in lemma \ref{MIPF} of the way a penny-farthing is connected with the surrounding ray category because one has to exclude the numbers 11,12,20,21 and 24.

Several times we use the following well-known result ( see e.g. \cite{trennen,Bongo} ). If in a ray category $P$ the composition of all arrows going through a point $x$ vanishes one can split this point into an emitter and a receiver to obtain a new quiver and a new ray category $P'$ with the obvious induced relations. Then $P$ is minimal representation-infinite iff $P'$ is so. When we apply this reduction to all points of the quotient by the square of the radical we simply say that we consider the separated quiver even though this quiver might still be connected.

Now we explain the simple strategy of the proof of the structure theorem thereby fixing some notations. For a non-deep contour $C=\{v,w\}$ we choose paths $v=v_{1}\ldots v_{n}$ and $w=w_{1}\ldots w_{m}$ from $x$ to $y$. Up to duality it suffices to consider the case that  $P(x,y)$ is a cyclic $P(y,y)$-module and  we choose a path $p=p_{1}p_{2}\ldots p_{r}$ such that $\rho=\vec{p}=\pi(y)$. We use the abbreviations $\alpha=\vec{v}_{1}$,$\beta=\vec{w}_{1}$,$\gamma=\vec{v}_{2}\ldots \vec{v}_{n}$, $\delta=\vec{w}_{2}\ldots \vec{w}_{m}$ and $P'=P/\rho\vec{v}$. Then the contour induces in $P$  the  diagram shown in figure 1.2. 
\vspace{0.5cm}

\setlength{\unitlength}{0.5cm}
\begin{picture}(15,6)
\put(8,5){\vector(4,1){4}}
\put(8,5){\vector(3,-1){3}}
\put(12,6){\vector(1,-1){1}}
\put(11,4){\vector(2,1){2}}
\put(8,2){\vector(4,1){4}}
\put(8,2){\vector(3,-1){3}}
\put(12,3){\vector(1,-1){1}}
\put(11,1){\vector(2,1){2}}
\put(13,5){\vector(0,-1){3}}
\put(8,-1){figure 1.2}
\put(10,6){$\gamma$}
\put(12.3,5.8){$\alpha$}
\put(10,3){$\gamma$}
\put(12.3,2.8){$\alpha$}
\put(9,4){$\delta$}
\put(9,1){$\delta$}
\put(12,4){$\beta$}
\put(12,1){$\beta$}
\put(13.5,3.5){$\rho$}
\end{picture}

\vspace{1cm}

The proof of theorem 1 is just a careful analysis of the fact that the obvious $\tilde{D}_{5}$ subdiagram or some diagrams deduced from it cannot be cleaving in $P'$ if the contour is mild. This analysis is rather straightforward and much easier than the proof of theorem 2 in \cite{gaps} because there the problem is not a local one and also deep contours have to be considered.

 As a first easy step we end this section by showing the uniqueness of the paths $v$ and $w$. In fact, also $p$ is uniquely determined if one admits decisive categories with 5 points, but we have no neeed for this result.

\begin{lemma}\label{wegeindeutig} We keep the notations introduced above and assume that $C$ is a mild contour such that $\vec{v}$ is transit. Then any path $u=u_{1}u_{2}\ldots u_{r}$  with $\vec{u}=\vec{v}$ is equal to $v$ or $w$.

\end{lemma}

Proof:  The first arrow $u_{1}$ cannot be different from $v_{1}$ and $w_{1}$ because otherwise $P'$ contains the first cleaving diagram from figure 1.3. So assume $u_{1}=v_{1}$ and choose $i$ maximal such that $u_{j}=v_{j}$ for all $j\leq i$ and set $\sigma = \vec{u}_{1}\ldots \vec{u}_{i}.$  None of $u$ or $v$ is a proper subpath of the other, but they are interlaced.  If $u$ and $v$ are different we get the second diagram in figure 1.3. This is cleaving since $w$ is not interlaced with $v$ and it lies in $P'$.

\setlength{\unitlength}{0.8cm}
\begin{picture}(15,3)

\put(7,1){\vector(1,0){1}}
\put(9,1){\vector(-1,0){1}}
\put(8,2){\vector(0,-1){1}}
\put(8,1){\vector(0,-1){1}}
\put(2.1,0.4){$\rho$}
\put(2.5,0.7){$\vec{u}_{1}$}

\put(1.4,1.3){$\alpha$}
\put(2.1,1.3){$\beta$}

\put(5,-1){figure 1.3}

\put(1,1){\vector(1,0){1}}
\put(3,1){\vector(-1,0){1}}
\put(2,2){\vector(0,-1){1}}
\put(2,1){\vector(0,-1){1}}
\put(8.1,0.4){$\rho$}
\put(8.4,0.7){$\sigma \vec{w}_{i+1}$}

\put(7.4,1.3){$\beta$}
\put(8.1,1.3){$\sigma \vec{u}_{i+1}$}

\end{picture}
\vspace{0.5cm}\newpage

\section{The proof of the structure theorem}
\subsection{The subdivision into three cases }

In the article on multiplicative bases the rough classification of the non-deep contours into three cases is based on Roiters transit-lemma. The next crucial result leads also to a trichotomy.

\begin{lemma}\label{factor}
Let $C=\{v,w\}$ be a mild contour as before such that $\vec{v}$ is transit. Then we are in one of the following situations:
\begin{enumerate}
 \item $\rho$ factors through  $\alpha$ or $\beta$.
\item  $\vec{v}$ is also cotransit and the generator $\sigma$ of the radical of $P(x,x)$ factors through one of the morphisms  $\vec{v}_{p}$ or $\vec{w}_{q}$.
\end{enumerate}
\end{lemma}

Proof: Assume to the contrary that none of the assertions is true. Then  we will find in $P'$ appropriate cleaving diagrams that are not representation-finite. 

First we have $\rho^{2}=0$ because otherwise the first diagram of figure 2.1 is cleaving. The four morphisms $\alpha$,$\beta$,$\rho\alpha$ and $\rho\beta$ in $P'$ all with the same codomain cannot induce a cleaving diagram of type $\tilde{D}_{4}$. Thus we are up to symmetry in $\alpha$ and $\beta$ in one of the cases  I:$\rho\alpha = \alpha\tau$ or II:$\rho\alpha = \beta\phi$. Here $\tau$ and $\phi$ are not identities because $\rho\alpha$ is not irreducible.

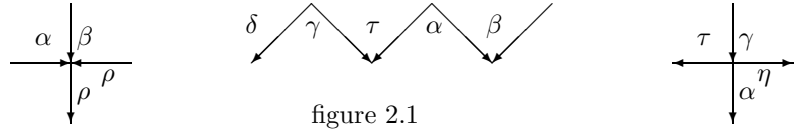
\setlength{\unitlength}{0.8cm}
\begin{picture}(15,3)

\put(13,1){\vector(0,-1){1}}
\put(13,1){\vector(-1,0){1}}
\put(13,1){\vector(1,0){1}}
\put(13,2){\vector(0,-1){1}}
\put(2.1,0.4){$\rho$}
\put(2.5,0.7){$\rho$}

\put(1.4,1.3){$\alpha$}
\put(2.1,1.3){$\beta$}

\put(6,0){figure 2.1}

\put(1,1){\vector(1,0){1}}
\put(3,1){\vector(-1,0){1}}
\put(2,2){\vector(0,-1){1}}
\put(2,1){\vector(0,-1){1}}
\put(13.1,0.4){$\alpha$}
\put(13.4,0.7){$\eta$}

\put(12.4,1.3){$\tau$}
\put(13.1,1.3){$\gamma$}\put(6,2){\vector(1,-1){1}}

\put(6,2){\vector(-1,-1){1}}\put(4.9,1.5){$\delta$}

\put(8,2){\vector(-1,-1){1}}\put(5.9,1.5){$\gamma$}

\put(8,2){\vector(1,-1){1}}\put(6.9,1.5){$\tau$}
\put(10,2){\vector(-1,-1){1}}\put(7.9,1.5){$\alpha$}
\put(8.9,1.5){$\beta$}

\end{picture}\vspace{0.5cm}

We analyze in detail the case $\rho\alpha=\alpha\tau$.
We claim that the middle diagram  in figure 2.1 is cleaving. Indeed, $\alpha$ and $\beta$ are irreducible and different. From $\tau= \xi\alpha$ one gets $\alpha\tau=\alpha\xi\alpha=\rho\alpha$ whence by cancellation  $\rho=\alpha\xi$. But $\rho$ does not factor. The relation $\alpha=\xi\tau$ implies $\alpha=\tau$ and $\rho=\alpha$, a contradiction. Next, $\tau=\gamma\xi$ gives $0\neq \rho\alpha\gamma= \alpha\tau\gamma=\alpha\gamma\xi\gamma$. Thus $\alpha\gamma$ is also cotransit and $\sigma=\xi\gamma$ factors. From $\gamma=\tau\xi$ we conclude $0\neq \rho\alpha\gamma=\rho\alpha\tau\xi=\rho^{2}\alpha\xi$ contradicting $\rho^{2}=0$. Finally none of $\delta$ or $\gamma$ factors through the other by lemma \ref{wegeindeutig}.

Now we look at  $\rho\beta$. For $\rho\beta= \alpha\psi$ the diagram $\tilde{A}(\beta,\psi,\tau,\alpha)$ is cleaving. To prove this only the possibilities $\psi=\xi\beta$, $\psi=\tau\xi$ and $\tau=\psi\xi$ remain to be excluded. But $\psi=\xi\beta$ implies $\alpha\psi=\alpha\xi\beta=\rho\beta$ whence $\rho=\alpha\xi$. The second gives $\alpha\psi=\alpha\tau\xi=\rho\beta=\rho\alpha\xi$ and therefore $\beta=\alpha\xi$. Similarly the third has $\alpha=\beta\xi$ as a consequence.

For $\rho\beta= \beta\psi$ the diagram $\tilde{A}(\beta,\psi,\delta,\gamma,\tau,\alpha)$ is cleaving. This time only  the possibilities $\psi=\xi\beta$, $\psi=\delta\xi$ and $\delta=\psi\xi$ have to be excluded. The first gives $\beta\psi=\beta\xi\beta=\rho\beta$ and $\rho= \beta\xi$ by cancellation. The second implies $\beta\psi\delta=\beta\delta\xi\delta=\rho\beta\delta$. Thus $\beta\delta$ is cotransit and we are in the second situation of the lemma. Finally the third possibility leads to $0\neq \rho\beta\delta=\rho\beta\psi\xi=\rho^{2}\beta\xi=0$.

We are reduced to the case that $\rho\beta$ does not factor through $\alpha$ or $\beta$ and  we consider the first representation-infinite category in figure 2.2 which is cleaving except for $\tau\gamma=\gamma\xi$. Then we look at the second diagram of figure 2.2.
\setlength{\unitlength}{0.8cm}
\begin{picture}(15,3)

\put(9,1){\vector(-1,0){1}}
\put(10,1){\vector(-1,0){1}}\put(10,1){\vector(1,0){1}}
\put(10,1){\vector(0,1){1}}

\put(12,1){\vector(1,0){1}}
\put(12,1){\vector(-1,0){1}}
\put(12,1){\vector(1,0){1}}
\put(14,1){\vector(-1,0){1}}
\put(15,1){\vector(-1,0){1}}
\put(10,1){\circle*{0.1}}\put(2,1){\circle*{0.1}}
\put(12,1){\circle*{0.1}}

\put(6,1){\vector(-1,0){1}}
\put(7,1){\vector(-1,0){1}}

\put(2,1){\vector(-1,0){1}}
\put(2,1){\vector(1,0){1}}
\put(4,1){\vector(1,0){1}}
\put(4,1){\vector(-1,0){1}}
\put(4,2){\vector(0,-1){1}}
\put(1.5,0.5){$\delta$}
\put(6.5,0.5){$\beta$}
\put(3.5,0.5){$\tau$}
\put(5.5,0.5){$\rho$}

\put(2.5,0.5){$\gamma$}
\put(4.2,1.5){$\gamma$}
\put(4.5,0.5){$\alpha$}

\put(8.5,0.5){$\delta$}
\put(9.5,0.5){$\xi$}
\put(10.5,0.5){$\gamma$}

\put(11.5,0.5){$\tau$}
\put(12.5,0.5){$\alpha$}
\put(13.5,0.5){$\rho$}
\put(14.5,0.5){$\beta$}
\put(10.3,1,5){$\delta$}
\put(6.5,-0.5){figure 2.2}

\end{picture}
\vspace{0.5cm}

If this is not cleaving we are in one of the following four situations. First $\delta=\eta\xi$ which implies $0\neq \rho\beta\delta=\rho\beta\eta\xi=\rho\alpha\gamma=\alpha\tau\gamma=\alpha\gamma\xi$ whence $\rho\beta\eta=\alpha\gamma$ and $\rho^{2}\neq0$. The second possibility $\gamma=\eta\xi$ leads to the same contradiction. The third possibility is $\delta\xi= \eta\delta$. This gives $\beta\delta\xi= \beta\eta\delta=\alpha\gamma\xi=\alpha\tau\gamma=\rho\alpha\gamma=\rho\beta\delta$ whence $\rho\beta=\beta\eta$ and we are in a case already dealt with. In the last case we have $\delta\xi=\eta\gamma$. We obtain $\beta\delta\xi=\beta\eta\gamma=\alpha\gamma\xi=\alpha\tau\gamma=\rho\alpha\gamma$ and finally $\rho\alpha=\beta\eta$. 

Thus we are in the situation that I and II are true at the same time and we show that this is impossible because the third diagram in figure 2.1 is cleaving.  Namely $0\neq \rho\alpha\gamma=\beta\eta\gamma =\alpha\tau\gamma$ shows that no composition in that diagram vanishes. The situation is symmetric in $\tau$ and $\eta$. So it only remains to exclude  $\tau=\xi\eta$. We obtain $\rho\alpha=\alpha\tau=\alpha\xi\eta=\beta\eta$ whence the contradiction $\beta=\alpha\xi$. 

The proof of the lemma is now complete if we are in case I. In case II the situation is not quite symmetric, but all arguments can be easily adopted. So we omit the details.

The resulting trichotomy reads as follows:

\begin{lemma} \label{cases}Let $C= \{v,w\}$ be a mild contour as before. Assume that $\vec{v}$ is transit and that $\rho=\alpha\omega$. Then $\rho$ does not factor through $\beta$ and we are in exactly one of the following three cases leading for $P(C)$ to a dumb-bell, a penny-farthing or  a diamond.
\begin{enumerate}
 \item (DB): $\omega$ is an identity and $\vec{v}_{2}=\beta$.
\item (PF): $\omega$ is an identity and $\vec{v}_{2}=\alpha$.
\item (D): $\omega$ is not an identity.
\end{enumerate}

\end{lemma}

Proof: At the end of the last proof we have seen that the cases I and II of that proof exclude each other. In particular, $\rho$ does not factor through $\beta$ and $\alpha$. If $\omega$ is an identity then we get $\alpha\vec{v}_{2}\neq 0$, $\alpha^{2}\neq 0$ and $\alpha\beta \neq 0$. Looking at the obvious diagram of type $\tilde{D}_{4}$ we conclude that $\vec{v}_{2}$ coincides with $\alpha$ or $\beta$.\newpage
\subsection{Dumb-bells} 

\begin{proposition}\label{PDB} In the case (DB) of lemma \ref{cases} $P(C)$ is a dumb-bell.
 
\end{proposition}

Proof: Let $v=v_{1}\ldots v_{n}$ and $w=w_{1}\ldots w_{m}$ be the corresponding uniquely determined paths from $x$ to $y$. By assumption we have $v_{2}=w_{1}$. Let $i\leq m$ be the maximal index with $w_{j}=v_{j+1}$ for all $j\leq i$. Set $\phi=\vec{w}_{2}\ldots\vec{w}_{i}=\vec{v}_{3}\ldots\vec{v}_{i+1}$,$\psi=\vec{v}_{n}\ldots\vec{v}_{i+2}$ and $\chi=\vec{w}_{i+1}\ldots\vec{w}_{m}$. Here $\chi$ is an identity for $i=m$ and $\psi$ is one for $i=n-1$. In any case we get $\alpha\beta\phi\psi=\vec{v}=\vec{w}=\beta\phi\chi$. So we see that $\chi$ is not an identity. If the same holds for $\psi$ we obtain a contradiction. Namely, $\psi$ and $\chi$ both belong to $P(x,z)$ and so they are comparable. But the uniquely determined first arrows of $\psi$ and $\chi$ are different because $v$ and $w$ are not interlaced, and the last arrows are different by the definition of $i$. 

Thus we have $\alpha\beta\phi=\beta\phi\chi$. Because $\vec{v}$ is transit one has $\chi=\pi(x)$. For $x=y$ we obtain $\chi=\alpha$ and $\beta\phi=\alpha^{i}$ for some $i$. But then $\{v,w\}$ is not a contour.

Thus $x$ is different from $y$ and $\beta\phi=\alpha^{i}\pi(x,y)$ with $i=0$ because otherwise $v$ and $w$ are interlaced. If $\chi=\chi'\pi(x,y)$ then we obtain $\alpha\pi(x,y)=\pi(x,y)\chi'\pi(x,y)$ whence the contradiction $\alpha=\pi(x,y)\chi'$. So $\chi$ does not factor through $\pi(x,y)$. It follows that for $P(y,x)\neq 0$ the diagram $\tilde{A}(\chi,\pi(x,y),\alpha,\pi(y,x))$ is cleaving in $P'$. Thus $P(y,x)=0$.

If $\chi$ factors properly as $\chi_{2}\chi_{1}$ the first diagram in figure 2.3 is cleaving. Namely $\chi^{2}\neq 0$ follows from $0\neq \alpha^{2}\beta\phi=\beta\phi\chi^{2}$ and $\beta\phi\chi=\alpha^{2}\xi$ implies that $v$ and $w$ are interlaced.

From $\phi\neq0$ we obtain the second cleaving diagram in figure 2.3. Here all possible factorizations contradict the fact that $v$ and $w$ are not interlaced.

Therefore the quiver of $Q$ is as it should be and the commutativity relation holds. The possible ray categories are then determined by the two nilpotence indices of the two loops. We use now the notations from the definition of a dumb-bell in figure 1.1. Up to duality we can assume $r\leq s$. Then we get $s\geq 3$ because the contour is non-deep and $s\leq 3$ because otherwise $P'$ contains the   category with number 14 from the BHV-list. For $r\geq 6$ one finds a category with number 21 in $P'$. The algebra with parameters $s=3$ and $r=5$ is representation-finite as follows easily from the criterion for finite representation type

\vspace{0.5cm}
\setlength{\unitlength}{0.8cm}
\begin{picture}(10,5).

\put(1,5){\vector(0,-1){1}}\put(0.4,4.5){$\chi_{1}$}\put(0.4,3.5){$\chi_{2}$}\put(0.4,2.5){$\chi_{1}$}\put(0.4,1.5){$\chi_{2}$}
\put(1,4){\vector(0,-1){1}}\put(2,3.5){$\beta\phi$} \put(3.5,3.5){$\alpha^{2}$}\put(8,3.5){$\beta\phi$}\put(8.5,4.5){$\beta$}
\put(9,5){\vector(0,-1){1}}\put(6.5,2.5){$\chi$}\put(6.5,1.5){$\chi$}\put(7.5,1.5){$\phi$}\put(9.5,3.5){$\alpha$}\put(9.5,4.5){$\alpha$}
\put(9,4){\vector(0,-1){1}}
\put(1,3){\vector(0,-1){1}}
\put(1,2){\vector(0,-1){1}}
\put(7,3){\vector(0,-1){1}}
\put(7,2){\vector(0,-1){1}}
\put(3,5){\vector(0,-1){2}}
\put(1,3){\vector(1,0){2}}

\put(7,2){\vector(1,0){1}}
\put(7,3){\vector(1,0){2}}

\put(8,4){\vector(1,0){1}}

\put(5,0){figure 2.3}\end{picture}\vspace{0.5cm}

\subsection{Penny-farthings}

\begin{proposition}\label{PPF}In the case (PF) of lemma \ref{cases} $P(C)$ is a penny-farthing.
 
\end{proposition}

Proof: Let $v=v_{1}\ldots v_{n}$ and $w=w_{1}\ldots w_{m}$ be the corresponding uniquely determined paths from $x$ to $y$. By assumption we have $\vec{v}_{2}=\alpha$. Let $i\geq 2$ be the maximal index with $\alpha=\vec{v}_{j}$ for all $j\leq i$ and assume $i<n$. Then we have $\vec{v}_{i+1}=\beta$ and the first diagram in figure 2.4 is in $P'$. Since it is not cleaving we have $\alpha^{2}=\beta\eta$ and $v$ and $w$ are interlaced. So $i=n$ and $\alpha^{n}=\beta\delta$. For $n\geq 3$ the second diagram in 2.4 is obviously cleaving whence $n=2$. The third diagram  in figure 2.4 enforces $\vec{w}_{m}\vec{w}_{1}=0$ because $0\neq \alpha^{3}=\alpha\vec{w}=\vec{w}\alpha$.
\setlength{\unitlength}{0.8cm}
\begin{picture}(10,5)

\put(1,4){\vector(0,-1){1}}\put(1.2,3.5){$\alpha$}\put(1.2,2.5){$\alpha$}\put(1.2,1.5){$\alpha$}
\put(4.5,3.5){$\alpha$}\put(4.5,2.5){$\alpha$}\put(5.5,1.5){$\alpha$}
\put(0.5,3.5){$\beta$}\put(0.5,2.5){$\beta$}\put(4.5,1.5){$\beta$}\put(5.5,3.5){$\delta$}
\put(1,3){\vector(0,-1){1}}\put(9.5,1.5){$\alpha$}\put(9.5,2.5){$\alpha$}
\put(1,2){\vector(0,-1){1}}\put(7.5,2.5){$\vec{w}_{1}$}\put(7.5,1.5){$\vec{w}_{m}$}

\put(5,4){\vector(0,-1){1}}
\put(5,3){\vector(0,-1){1}}
\put(5,2){\vector(0,-1){1}}

\put(0,3){\vector(1,0){1}}
\put(0,2){\vector(1,0){1}}

\put(5,3){\vector(1,0){1}}
\put(4,2){\vector(1,0){1}}

\put(9,3){\vector(0,-1){1}}
\put(9,2){\vector(0,-1){1}}

\put(8,3){\vector(1,-1){1}}
\put(9,2){\vector(-1,-1){1}}
\put(5,0){figure 2.4}\end{picture}

For the rest of the proof we use the notations of figure 1.1. If $x_{i}=x_{0}$ for some $i\neq 0$ we have $\alpha_{i}\ldots \alpha_{1}=\rho^{r}$ and $v$ and $w$ are interlaced. Next assume $x_{i}=x_{j}=z$ for some $i$,$j$ with $0<i<j\leq n-1$. Up to duality $P(z,x)$ is transit and we get $\alpha_{n-1}\ldots\alpha_{j} \alpha_{j-1}\ldots\alpha_{i}= \rho^{r}\alpha_{n-1}\ldots\alpha_{j}$ whence the path $w$ is not uniquely determined as it should be. The same argument applies if there are any additional arrows within $Q$. So the points and the arrows of $Q$ are  those of a penny-farthing. Furthermore we have the wanted commutativity relation and the zero relation  $\alpha_{n-1}\alpha_{1}=0$. Then the only possible relations for a ray category are as given in the definition of a penny-farthing.

\subsection{Diamonds}

 The first diagram in figure 2.5 cannot be cleaving. Therefore we are in one of the  cases $\gamma=\omega\xi$,
 $\omega\beta=\gamma\eta$ 
or $\omega\alpha=\gamma\eta$ with a non-identity $\eta$.\vspace{0.5cm}

\setlength{\unitlength}{0.8cm}
\begin{picture}(20,5)

\put(1,5.5){\vector(2,-1){1}}\put(1.2,5){$\alpha$}
\put(1,4.5){\vector(2,1){1}}\put(1.2,4){$\beta$}
\put(2,5){\vector(-1,-2){1}}\put(1.6,3.5){$\omega$}

\put(0,2.5){\vector(2,1){1}}
\put(1,3){\vector(2,-1){1}}\put(1.2,2.5){$\alpha$}\put(0.6,2.5){$\gamma$}

\put(5,3){\vector(1,1){1}}\put(4.8,3.5){$\vec{v}_{n}$}
\put(6,4){\vector(1,0){1}}\put(4.8,2.4){$\vec{w}_{m}$}
\put(7,4){\vector(1,0){1}}
\multiput(8.5,4)(0.5,0){7}{\circle*{0.05}}
\multiput(8.5,2)(0.5,0){7}{\circle*{0.05}}
\put(12,4){\vector(1,-1){1}}
\put(13,3){\vector(1,0){1}}\put(13.5,2.5){$\omega$}
\put(5,3){\vector(1,-1){1}}\put(11.8,3.5){$\vec{v}_{1}$}
\put(6,2){\vector(1,0){1}}\put(11.8,2.4){$\vec{w}_{1}$}
\put(7,2){\vector(1,0){1}}

\put(12,2){\vector(1,1){1}}
\put(5,1){figure 2.5}\end{picture}

\begin{lemma}\label{LDI}
 The case $\gamma=\omega\xi$ is impossible.
\end{lemma}

Proof: For  $\gamma=\omega$ we get $0\neq \alpha\gamma\alpha\gamma = \beta\delta\alpha\gamma=\alpha\gamma\beta\delta=\beta\delta\beta\delta$ whence $\alpha$,$\beta$,$\gamma$ and $\delta$ induce an obvious cleaving diagram of type $\tilde{D}_{4}$ in $P'$. Thus $\xi$ is not an identity.

By definition we have $\omega=\vec{p}_{2}\ldots \vec{p}_{r}=\vec{v}_{2}\ldots \vec{v}_{r}$ and $\delta=\vec{w}_{2}\ldots \vec{w}_{m}$. Also the second diagram of figure 2.5 is cleaving in $P$.
Since $P'$ contains not the obvious  subdiagrams of types $\tilde{E}_{7}$ resp. $\tilde{E}_{6}$  we see that $\omega$ and $\delta$ are irreducible. Thus we have $r=m=2<n$. We look at the $\tilde{D}_{4}$-diagram with $\alpha$,$\beta$,$\omega$ and $\vec{v}_{3}$ as irreducible morphisms and conclude $\vec{v}_{3}=\beta $ or $\vec{v}_{3}=\alpha$. 

So we have in the first case $\gamma= \omega\beta\zeta$. For $\zeta \neq id$ one has by axiom e) of a ray category that $\zeta=\phi\delta$ - leading to the contradiction $0\neq \beta\delta=\alpha\omega\beta\phi\delta$ - or $\zeta=\delta\phi$ implying $0\neq \beta\delta=\alpha\omega\beta\delta\phi$. Therefore we have $\beta\delta= \alpha\omega\beta$ and $v$ is bitransit. We are in the situation dual to the case (DB) in lemma  \ref{cases} and we find the contradiction that $\omega$ is an identity.

In the remaining case $\gamma= \omega\alpha\zeta$ we assume first $n=3$ i.e. $\zeta= id$. Then  the first diagram of figure 2.6 is cleaving because $0\neq \alpha\omega\alpha\omega\alpha=\beta\delta\omega\alpha$. For $n\geq 4$ and $\vec{v}_{4}\neq \omega$ the second diagram of figure 2.6 is cleaving. Here  $\alpha\omega=\beta\phi$ contradicts lemma 1 and $\alpha\vec{v}_{4}=\beta\phi$ leads again to the case $\gamma=\omega\beta\zeta$ already dealt with. Thus we can assume $\vec{v_{4}}=\omega$. For $n=4$ the third diagram in figure 2.6 is cleaving. For $n= 5$ we distinguish  the cases $\vec{v}_{5}\neq \alpha$ and $\vec{v}_{5}= \alpha$. In the first case, the fourth diagram in figure 2.6 is cleaving, in the second case the first diagram in figure 2.5.

\vspace{0.5cm}
\setlength{\unitlength}{0.8cm}
\begin{picture}(20,2)

\put(0,1){\vector(1,0){1}}\put(1.5,0.5){$\alpha$}
\put(1,1){\vector(1,0){1}}\put(0.5,0.5){$\omega$}
\put(2,1){\vector(1,0){1}}\put(2.5,0.5){$\omega$}
\put(1,1){\vector(0,1){1}}\put(2.5,1.5){$\beta$}
\put(2,2){\vector(0,-1){1}}\put(0.5,1.5){$\delta$}

\put(4,1){\vector(1,0){1}}\put(4.5,0.5){$\omega$}
\put(5,1){\vector(1,0){1}}\put(5.5,0.5){$\alpha$}
\put(6,1){\vector(1,0){1}}\put(6.5,0.5){$\omega$}
\put(5,2){\vector(0,-1){1}}\put(6.5,1.5){$\beta$}
\put(6,2){\vector(0,-1){1}}\put(4.5,1.5){$\vec{v}_{4}$}

\put(8,1){\vector(1,0){1}}\put(8.5,0.5){$\alpha$}
\put(9,1){\vector(1,0){1}}\put(9.3,0.5){$\alpha\omega$}
\put(10,1){\vector(1,0){1}}\put(10.5,0.5){$\omega$}
\put(9,1){\vector(0,1){1}}\put(10.5,1.5){$\beta$}
\put(10,2){\vector(0,-1){1}}\put(8.5,1.5){$\delta$}

\put(12,1){\vector(1,0){1}}\put(12.5,0.5){$\alpha$}
\put(13,1){\vector(1,0){1}}\put(13.3,0.5){$\alpha\omega$}
\put(14,1){\vector(1,0){1}}\put(14.5,0.5){$\omega$}
\put(13,2){\vector(0,-1){1}}\put(14.5,1.5){$\beta$}
\put(14,2){\vector(0,-1){1}}\put(12.5,1.5){$\vec{v}_{5}$}

\put(6,-0.5){figure 2.6}\end{picture}\vspace{0.5cm}

\begin{lemma}\label{LDII}
 The case 
 $\omega\beta=\gamma\eta$ leads to the impossible case $\gamma=\omega\xi$ or to the case
 $\omega\alpha=\gamma\phi$.
\end{lemma}

Proof: We can assume that $\eta$ is not an identity and we claim that the diagram in figure 2.7 is cleaving. 

\vspace{0.5cm}
\setlength{\unitlength}{0.8cm}
\begin{picture}(10,2)

\put(2,1){\vector(-1,0){1}}
\put(3,1){\vector(-1,0){1}}
\put(3,1){\vector(1,0){1}}
\put(3,2){\vector(0,-1){1}}
\put(3,2){\vector(1,0){1}}
\put(4,2){\vector(0,-1){1}}
\put(5,2){\vector(-1,0){1}}
\put(6,2){\vector(-1,0){1}}
\put(5.5,2.1){$\alpha$}\put(4.5,2.1){$\omega$}\put(3.5,2.1){$\gamma$}
\put(1.5,0.5){$\delta$}\put(2.5,0.5){$\eta$}\put(3.5,0.5){$\beta$}
\put(2.5,1.5){$\delta$}\put(4.2,1.5){$\alpha$}

\put(5,0){figure 2.7}\end{picture}

\vspace{0.5cm}
First we get from  $0\neq\alpha\omega\beta\delta=\alpha\gamma\eta\delta=\beta\delta\eta\delta=\alpha\omega\alpha\gamma$ that $\delta\eta\delta\neq 0 \neq \alpha\omega\alpha$. We have to exclude many possible relations. To begin with $\omega\alpha= \gamma\theta$ leads us to the third case and $\gamma=\omega\theta$ to the first. 

So assume $\alpha\omega\alpha=\beta\theta$. Then $\tilde{A}(\delta,\gamma,\omega\alpha,\theta)$ is cleaving. For $\theta=\xi\omega\alpha$ implies $\beta\theta=\beta\xi\omega\alpha=\alpha\omega\alpha$ and thus $\alpha=\beta\xi$. In the same way $\omega\alpha=\xi\theta$ leads to $\beta=\alpha\xi$. The relation $\delta=\theta\xi$ gives $\beta\delta=\beta\theta\xi=\alpha\omega\alpha\xi=\alpha\gamma$ and we end up in the first case.  Similarly $\theta=\delta\xi$ implies $\gamma\xi=\omega\alpha$. 

The remaining part of the proof is almost dual, but we give full details. From $\delta\eta=\xi\beta$ we obtain $\beta\delta\eta\delta=\beta\xi\beta\delta=\alpha\gamma\eta\delta=\alpha\omega\beta\delta\neq 0$, whence $\rho=\alpha\omega=\beta\xi$ factors also through $\beta$. From $\beta=\xi\eta$ one gets $\gamma=\xi\omega$.

Finally we have to exclude $\delta\eta\delta=\theta\gamma$. Here $\theta$ is not an identity because of lemma \ref{wegeindeutig}. We claim that $\tilde{A}(\beta,\delta\eta,\theta,\alpha)$ is cleaving. Here $\theta=\alpha\xi$ implies $\theta\gamma=\alpha\xi\gamma=\delta\eta\delta=\xi\beta\delta$ whence $\delta\eta=\xi\beta$ which was just shown to be impossible. From $\theta=\delta\eta\xi$ one finds $\delta=\xi\gamma$ and similarly one gets $\gamma=\xi\delta$ from $\delta\eta=\theta\xi$. The proof of the second lemma is complete.

\begin{proposition}\label{PD}
 In the case $\omega\alpha=\gamma\eta$ with $\eta\neq id$ the contour is a diamond. Furthermore, $\beta$ is the only arrow $\epsilon$ with $\epsilon\delta\neq 0$.
\end{proposition}

Proof: We look at the diagram shown in figure 2.8 and we show first that some morphisms therein do not factor through each other. 

From $\delta\eta=\theta\alpha$ one gets $\beta\delta\eta=\beta\theta\alpha=\alpha\omega\alpha$. Thus $\rho=\alpha\omega=\beta\theta$ factors also through $\beta$. If $\alpha= \theta\eta$ then $\alpha=\eta$ and $\omega=\gamma$. Hence $\alpha$,$\beta$,$\gamma$ and $\delta$ induce a $\tilde{D}_{4}$-diagram in $P'$.

\setlength{\unitlength}{0.8cm}
\begin{picture}(10,4)

\put(3,2){\vector(1,1){1}}
\put(3,2){\vector(1,-1){1}}
\put(4,1){\vector(1,1){1}}
\put(4,3){\vector(1,-1){1}}
\put(4,3){\vector(2,0){2}}
\put(5,2){\vector(2,0){2}}
\put(6,3){\vector(1,1){1}}
\put(6,3){\vector(1,-1){1}}

\put(7,2){\vector(1,1){1}}
\put(7,4){\vector(1,-1){1}}
\put(3,2.3){$\gamma$}\put(4,2.3){$\alpha$}\put(6,2.3){$\gamma$}\put(7.8,2.3){$\alpha$}
\put(3,1.3){$\delta$}\put(4,1.3){$\beta$}\put(6,3.3){$\delta$}\put(7.8,3.3){$\beta$}
\put(5,3.3){$\eta$}\put(5.9,1.5){$\omega$}

\put(5,0){figure 2.8}\end{picture}\vspace{0.5cm}

If we have $\delta\eta\gamma=\theta\delta$ then $\theta$ is not an identity. We find in $P'$ the cleaving diagram $\tilde{A}(\alpha,\delta\eta,\theta,\beta)$. Namely $\delta\eta=\theta\xi$ gives $\beta\delta\eta\gamma=\beta\theta\xi\gamma=\beta\theta\delta$ and so $\delta=\xi\gamma$. Similarly, $\theta=\delta\eta\xi$ implies $\xi\delta=\gamma$. From $\theta= \xi\beta$ we obtain $\theta\delta= \xi\beta\delta=\delta\eta\gamma=\xi\alpha\gamma$ and $\delta\eta=\xi\alpha$ which is impossible as already seen. The case $\beta=\xi\theta$ leads to $\beta=\theta$ which was considered just before.

Finally set $\delta'=\vec{w}_{1}\ldots \vec{w}_{i}$ and $\delta''=\vec{w}_{i+1}\ldots \vec{w}_{m}$ for some $1\leq i \leq m-1$. Suppose $\omega\delta'=\gamma\theta.$ Then $\tilde{A}(\eta,\alpha,\delta',\theta)$ is cleaving in $P'$. Indeed $\theta=\eta\xi$ gives $\gamma\theta\delta''=\gamma\eta\xi\delta''=\omega\alpha\xi\delta''=\omega\delta'\delta''$ and therefore $\alpha\xi=\delta'$ contradicting lemma 1.
In the same vein $\eta=\theta\xi$ implies $\alpha=\delta'\xi$. From $\eta=\xi\alpha$ one comes to the already excluded case $\delta\eta=\delta\xi\alpha$ and $\alpha=\xi\eta$ gives $\alpha=\eta$. The cases $\alpha=\delta'\xi$ and $\delta'=\alpha\xi$ contradict lemma 1.
So assume $\delta'=\xi\theta$. This implies $\omega\delta'=\omega\xi\theta=\gamma\theta$ whence $\gamma=\omega\xi$. Since $\rho=\alpha\omega$ is a generator of the radical of $P(y,y)$ we obtain $\alpha\gamma=\alpha\omega\xi=(\alpha\omega)^{i}$ for some $i\geq 1$. We infer $\xi= id$ and $\gamma=\omega$. But now $\alpha$,$\beta$,$\gamma$ and $\delta$ produce a cleaving diagram of type $\tilde{D}_{4}$. Finally look at $\theta= \xi\delta'$. We get $\gamma\theta=\gamma\xi\delta'=\omega\delta'$ and $\gamma\xi=\omega$. From $\alpha\gamma=(\alpha\omega)^{i}=(\alpha\gamma\xi)^{i}$ we infer again $\gamma=\omega$ and we are back in the last case.

After these lengthy considerations we start with the proof. First we claim that $\gamma$ is irreducible in $P$. If not, the first diagram in figure 2.9 is cleaving in $P'$ by the above considerations. Similarly, if $\delta$ is not irreducible the second diagram in figure 2.9 is cleaving. Here $\xi=\vec{w}_{2}\ldots \vec{w}_{m-1}$.

\vspace{0.5cm}
\setlength{\unitlength}{0.8cm}
\begin{picture}(20,3)

\put(0,2){\vector(1,0){1}}\put(0.5,2.2){$\vec{v}_{2}$}
\put(1,2){\vector(1,0){1}}\put(1.5,2.2){$\eta$}
\put(2,2){\vector(1,0){1}}\put(2.5,2.2){$\delta$}
\put(0,1){\vector(1,0){1}}\put(0.5,0.5){$\beta$}
\put(1,1){\vector(1,0){1}}\put(1.5,0.5){$\omega$}
\put(1,2){\vector(0,-1){1}}\put(0.5,1.5){$\alpha$}\put(2.2,1.5){$\gamma$}
\put(2,2){\vector(0,-1){1}}

\put(6.5,2.2){$\eta$}\put(8.5,2.2){$\xi$}\put(4.5,0.5){$\xi$}
\put(7.3,2.2){$\vec{w}_{m}$}
\put(5.5,0.5){$\beta$}
\put(6.5,0.5){$\omega$}
\put(5.5,1.5){$\alpha$}
\put(7.2,1.5){$\gamma$}

\put(6,2){\vector(1,0){1}}
\put(6,1){\vector(1,0){1}}
\put(7,2){\vector(1,0){1}}
\put(8,2){\vector(1,0){1}}
\put(4,1){\vector(1,0){1}}
\put(5,1){\vector(1,0){1}}
\put(6,2){\vector(0,-1){1}}
\put(7,2){\vector(0,-1){1}}

\put(10,3){\vector(1,0){1}}\put(10.5,3.2){$\beta$}
\put(11,3){\vector(1,0){1}}\put(11.5,3.2){$\lambda$}
\put(12,3){\vector(1,0){1}}\put(12.5,3.2){$\kappa$}
\put(13,3){\vector(1,0){1}}\put(13.5,3.2){$\delta$}

\put(11,2){\vector(1,1){1}}\put(11,2.5){$\gamma$}
\put(11,2){\vector(1,-1){1}}\put(12.6,2.5){$\alpha$}
\put(12,3){\vector(1,-1){1}}\put(11,1.2){$\delta$}
\put(12,1){\vector(1,1){1}}\put(12.6,1.2){$\beta$}
\put(12,1){\vector(1,0){1}}\put(12.5,0.5){$\epsilon$}\put(10,0.5){$\kappa\lambda=0$}

\put(5,-0.5){figure 2.9}\end{picture}\vspace{0.5cm}
Now $\rho\alpha\gamma=\alpha\gamma\eta\gamma$ shows that $\vec{v}$ is also cotransit and that $\eta\gamma$ is a generator of the radical of $P(x,x)$. Thus the contour $\{v,w\}$ is self - dual. It follows that the diagram in  figure 2.8 is cleaving. To see this one checks first with the aid of the above considerations that the cleaving conditions are satisfied in all sources. By duality the same holds for all sinks. The relation $\kappa\lambda =0$ follows by looking at the $\tilde{D}_{4}$-diagram induced by $\gamma,\alpha,\lambda$ and $\kappa$.

Next, $x=y$ implies $\alpha\gamma=(\alpha\omega)^{i}$ because $\rho$ is a generator of the radical. Since $\gamma$ is irreducible we get $\gamma=\omega$. This induces the familiar cleaving diagram induced by $\alpha$,$\beta$,$\gamma$  and $\delta$. The starting points $z$ and $t$ of $\alpha$ and $\beta$ are different from $y$ because otherwise $\rho$ would be a loup. Dually, $x$ is different from $z$ and $t$. Since $\alpha$ and $\beta$ are different so are $z$ and $t$. Now we also use the notations from figure 1.1, i.e. we set $\omega=\lambda$ and $\kappa=\eta$.

Another arrow $\epsilon$ in $P$ with $\epsilon\delta\neq 0$ gives the third diagram from figure 2.9 in $P'$. This is cleaving because $\phi\delta= \xi\gamma$ gives rise to $\tilde{A}(\alpha,\xi,\phi,\beta)$. We claim that $\pi=\pi(z)=\gamma\kappa$. Otherwise we have $\pi^{i}=\gamma\kappa$ for some $i\geq 2$ and one sees easily that there is a cleaving diagram of type $\tilde{D}_{4}$ with arrows $\gamma,\pi, \alpha$ and $\pi$ again. It follows that $\gamma$ is bitransit. By duality the same holds for $\alpha$.

Next we show that $\kappa$ is irreducible in the full subcategory $Q$ supported by $x,y,z,t$. Namely,  if $\kappa$ factors through $y$ we have  $\kappa=\xi\alpha$ ( $\alpha$ is bitransit ) and there is again a $\tilde{D}_{4}$-diagram with arrows $\gamma,\delta, \xi\alpha$ and $\xi\beta$ cleaving. A factorization $\kappa=\psi\phi$ through $t$ implies the contradiction $0\neq \phi\gamma=\delta\xi= \delta(\kappa\gamma)^{i}=\delta\psi\phi(\gamma\kappa)^{i-1}\gamma$ or $0\neq \phi\gamma=\xi\delta=\zeta\beta\delta=\zeta\alpha\gamma$ whence $\phi=\zeta\alpha$ and $\kappa$ factors through $y$ which is impossible. Thus there is an arrow $\kappa'$ from $z$ to $x$ and from $\kappa=\kappa'(\gamma\kappa)^{i}$ or $\kappa=(\kappa\gamma)^{i}\kappa'$ one concludes $\kappa=\kappa'$. By duality, $\lambda$ is also irreducible in $Q$.

We still have to prove that there are no more arrows in $Q$. Since $\alpha$ is an arrow, the relation $\alpha\gamma=\pi(x,y)(\kappa\gamma)^{i}$ implies $i=0$, i.e. $\alpha\gamma=\pi(x,y)$ and there is no arrow from $x$ to $y$. Let $\xi:y \rightarrow x$ be an arrow in $Q$. Then $\tilde{A}(\alpha\gamma,\kappa\gamma,\xi,\alpha\lambda)$ is cleaving. For $\kappa\gamma=\xi\epsilon$ gives $\kappa\gamma=\xi(\alpha\lambda)^{i}\alpha\gamma$ contradicting the irreducibility of $\kappa$. Similarly, $\alpha\lambda=\epsilon\xi$ implies that $\lambda$ is not irreducible. Next, an arrow $\xi:t \rightarrow x$ induces the cleaving diagram $\tilde{A}(\beta,\xi,\kappa,\alpha)$ and a loop in $x$ is impossible because of $\pi(x)=\kappa\gamma$. Thus $\delta,\gamma$ and $\kappa$ are the only arrows starting or ending in $x$. Looking at the separated quiver one gets that there is at most one more arrow and this flies from  $z$ to $t$ or  from $t$ to $z$ or from $t$ to $t$. If there is an arrow $\xi:z \rightarrow t$ each nilpotent endomorphism of $t$ is a power of $\xi\lambda\beta$ whence $P(z,t)$ is cyclic over $P(z,z)$. We obtain $\delta\kappa=\xi(\gamma\kappa)^{i}$ and the contradiction $\delta=\xi(\gamma\kappa)^{i-1}\gamma$. Dually there is no arrow from $t$ to $z$. Finally, a loop $\xi$ at $t$ induces the cleaving diagram $\tilde{A}(\xi,\beta,\alpha,\delta\kappa)$. Here $\delta\kappa=\xi\epsilon=\xi\delta\kappa$ is impossible.

\section{The proof of the disjointness theorem}

\subsection{Some preparations}

The proof of theorem 2 is much easier than that of theorem 1 and many arguments are the same as in \cite{BGRS}. But for the convenience of the reader and also because some of our statements are more general we give full details.
We start with three easy lemmata that will also be used in the next section.
\begin{lemma}\label{LDB}
 Let $C$ be a mild dumb-bell. Then $\mu$ and $\rho$ are the only arrows ending in $y$.
\end{lemma}

Proof: If not let $\nu:z \rightarrow y$ be an additional arrow. Then we have $\rho \nu =0$ because we get otherwise an obvious cleaving diagram of type $\tilde{D}_ {4}$. Consider the first diagram in figure 3.1. If this is not cleaving we have $\rho \mu = \nu\xi$ or $\rho^{2}=\nu\xi$ for some $\xi$. The first case implies the contradiction $0\neq \rho^{2}\mu= \rho\nu\xi= 0$. In the second case we write $\xi=\chi\eta$ for some arrow $\eta$ and we get $0\neq \rho^{2}\mu = \nu\chi\eta\mu$. Then the second diagram in figure 3.1 is cleaving unless $\eta\mu= \zeta\lambda$. But then $\tilde{A}(\rho,\eta,\zeta,\mu)$ is cleaving.

\vspace{0.5cm}
\setlength{\unitlength}{0.8cm}
\begin{picture}(15,3)

\put(4,2){\vector(-1,0){1}}
\put(2,2){\vector(1,0){1}}
\put(4,1){\vector(-1,0){1}}
\put(2,1){\vector(-1,0){1}}
\put(2,1){\vector(1,0){1}}
\put(2,2){\vector(0,-1){1}}
\put(3,2){\vector(0,-1){1}}
\put(3.5,0.5){$\nu$}
\put(2.5,0.5){$\mu$}
\put(1.5,0.5){$\lambda$}
\put(2.5,2.2){$\mu$}
\put(3.5,2.2){$\rho$}

\put(1.5,1.5){$\lambda$}
\put(2.5,1.5){$\rho$}

\put(9,2){\vector(1,0){1}}
\put(8,2){\vector(1,0){1}}
\put(10,1){\vector(-1,0){1}}
\put(8,1){\vector(-1,0){1}}
\put(8,1){\vector(1,0){1}}
\put(8,2){\vector(0,-1){1}}
\put(9,2){\vector(0,-1){1}}
\put(9.5,0.5){$\nu$}
\put(8.5,0.5){$\mu$}
\put(7.5,0.5){$\lambda$}
\put(8.5,2.2){$\mu$}
\put(9.5,2.2){$\eta$}

\put(7.5,1.5){$\lambda$}
\put(8.5,1.5){$\rho$}
\put(5,0){figure 3.1}\end{picture}\vspace{0.5cm}

\begin{lemma}\label{LPF}
 Let $C$ be a mild penny-farthing with $C(x_{0},y)\neq 0$ for some $y\notin C$. Then $n=2$.
\end{lemma}

Proof: Suppose $n>2$. First , let $\epsilon:x_{0} \rightarrow y$ be an additional arrow. Then we obtain $\epsilon\rho=0$ from an obvious $\tilde{D}_{4}$-diagram and  the first diagram in figure 3.2 is cleaving. Here $\alpha_{3}'=\alpha_{n}\ldots \alpha_{3}$. For $\xi\epsilon= \alpha_{2}\alpha_{1}\rho$ leads to the contradiction $0=\alpha_{3}'\xi\epsilon\rho=\alpha_{n}\ldots \alpha_{1}\rho=\rho^{3}\neq 0$ and $\xi\epsilon= \rho^{2}$  to $0=\xi\epsilon\rho=\rho^{3}\neq 0$. So from now on we assume that $\alpha_{1}$ and $\rho$ are the only arrows starting at $x_{0}$.

Choose $i$ minimal such that there is an arrow $\epsilon:x_{i} \rightarrow y$ with $\epsilon\alpha_{i}\ldots \alpha_{1}\neq 0$. Define $\alpha'= \alpha_{i}\ldots \alpha_{1}$ and $\alpha''= \alpha_{n}\ldots \alpha_{i+1}$. We consider the case $i>1$ and set $\alpha_{2}'=\alpha_{i}\ldots \alpha_{2}$. The case $i=1$ is analogous. We claim that
 the second diagram in figure 3.2 is cleaving. Namely $\alpha''=\xi\epsilon$ contradicts the uniqueness lemma. From $0 \neq \epsilon\alpha'=\xi \rho$ we obtain $\xi=\xi'\alpha'$ because of the minimal choice of $i$. Thus we get $\xi'=\epsilon\eta$ implying $0\neq \epsilon\alpha'=\epsilon\eta\alpha'\rho$ or $\xi'=\eta\epsilon$ implying $0\neq \epsilon\alpha'=\eta\epsilon\alpha'\rho$.

\vspace{0.5cm}
\setlength{\unitlength}{0.8cm}
\begin{picture}(15,4)

\put(2,3){\vector(0,1){1}}
\put(2,3){\vector(1,-1){1}}
\put(2,3){\vector(-1,-1){1}}
\put(1,2){\vector(1,-1){1}}
\put(3,2){\vector(-1,-1){1}}
\put(3,2){\vector(1,-1){1}}
\put(4,1){\vector(1,0){1}}

\put(4,3){\vector(-1,-1){1}}
\put(5,3){\vector(-1,0){1}}
\put(1.5,3.5){$\epsilon$}
\put(0.8,2.5){$\alpha_{1}$}
\put(0.4,1.5){$\alpha_{3}'\alpha_{2}$}
\put(2.8,2.5){$\rho$}
\put(2.8,1.5){$\rho$}
\put(3.8,1.5){$\alpha_{1}$}
\put(3.8,2.5){$\alpha_{3}'$}
\put(4.5,3.2){$\alpha_{2}$}
\put(4.5,0.6){$\alpha_{2}$}
\put(4.4,2){$\alpha_{1}\alpha_{3}'=0$}

\put(8,2){\vector(-1,-1){1}}
\put(9,3){\vector(1,-1){1}}
\put(9,3){\vector(-1,-1){1}}
\put(8,2){\vector(1,-1){1}}
\put(10,2){\vector(-1,-1){1}}
\put(10,2){\vector(1,-1){1}}
\put(11,1){\vector(1,0){1}}

\put(11,3){\vector(-1,-1){1}}
\put(12,3){\vector(-1,0){1}}

\put(7,1.5){$\epsilon$}
\put(7.8,2.5){$\alpha'$}
\put(7.8,1.2){$\alpha''$}
\put(9.8,2.5){$\rho$}
\put(9.8,1.5){$\rho$}
\put(10.8,1.5){$\alpha_{1}$}
\put(10.8,2.5){$\alpha''$}
\put(11.5,3.2){$\alpha_{2}'$}
\put(11.5,0.6){$\alpha_{2}'$}
\put(11.4,2){$\alpha_{1}\alpha''=0$}
\put(6.5,0){figure 3.2}\end{picture}\vspace{0.5cm}

\begin{lemma}\label{LD}
 In a mild diamond, $\beta$ is the only arrow whose composition with $\delta$ does not vanish and $\gamma$ and $\delta$ are the only arrows starting at $x$.
\end{lemma}
 
Proof: The first assertion has already been shown in proposition \ref{PD}. Let $\epsilon:x \rightarrow a$ be an additional arrow. If the diagram of figure 3.3 is not cleaving we have $\xi\epsilon=\alpha\gamma$ contradicting the uniqueness-lemma or $\xi\epsilon=\delta\kappa\gamma$ whence $\beta\xi\epsilon=\beta\delta\kappa\gamma\neq 0$ for some appropriate $\xi$. By the dual of the first assertion of the lemma we get $\xi=\delta\zeta$. But now the $\tilde{D}_{4}$-diagram with morphisms $\kappa,\zeta,\gamma,\delta$ is cleaving because the factorization of $\kappa$ by $\zeta$ resp. of $\zeta$ by $\kappa$ contradicts the irreducibility of $\epsilon$ resp. $\gamma$.

\vspace{0.5cm}
\setlength{\unitlength}{0.7cm}
\begin{picture}(15,4)

\put(2,3){\vector(0,1){1}}
\put(2,3){\vector(1,-1){1}}
\put(2,3){\vector(-1,-1){1}}
\put(1,2){\vector(1,-1){1}}
\put(3,2){\vector(-1,-1){1}}
\put(3,2){\vector(1,-1){1}}
\put(4,1){\vector(1,0){1}}

\put(4,3){\vector(-1,-1){1}}
\put(5,3){\vector(-1,0){1}}
\put(1.5,3.5){$\epsilon$}
\put(0.8,2.5){$\delta$}
\put(0.8,1.5){$\beta$}
\put(2.8,2.5){$\gamma$}
\put(2.8,1.5){$\alpha$}
\put(3.8,1.5){$\kappa$}
\put(3.8,2.5){$\lambda$}
\put(4.5,3.2){$\beta$}
\put(4.5,0.6){$\delta$}
\put(4.4,2){$\kappa\lambda=0$}

\put(6.5,0){figure 3.3}\end{picture}\vspace{0.5cm}

\subsection{The proof of theorem 2}

Now we prove the parts a) and b) of  theorem 2 by showing that a common arrow or a common point of two different mild contours always leads to a contradiction to the assumptions of a) or b). 

First, let $C$ be a dumb-bell from $x$ to $y$. Up to duality we can suppose that $y$  occurs in another mild contour $C'$. If $C'$ is another dumb-bell we have $x'=y$ by \ref{LDB}. Then the diagram in figure 3.4 is cleaving.

\vspace{0.5cm}
\setlength{\unitlength}{0.8cm}
\begin{picture}(15,3)

\put(2,1){\vector(-1,0){1}}
\put(3,1){\vector(-1,0){1}}
\put(3,1){\vector(1,0){1}}
\put(3,1){\circle*{0.1}}
\put(5,1){\vector(-1,0){1}}
\put(5,1){\vector(1,0){1}}
\put(5,2){\vector(0,-1){1}}
\put(5,2){\vector(1,0){1}}

\put(6,2){\vector(0,-1){1}}
\put(7,2){\vector(-1,0){1}}
\put(1.5,0.5){$\lambda$}
\put(2.5,0.5){$\lambda$}
\put(3.5,0.5){$\mu$}
\put(4.5,0.5){$\rho$}
\put(5.5,0.5){$\mu'$}
\put(5.5,2.5){$\mu'$}
\put(6.5,2.5){$\rho'$}
\put(6.5,1.5){$\rho'$}
\put(4.5,1.5){$\rho$}

\put(6.5,0){figure 3.4}
\end{picture}\vspace{0.5cm}

If $C'$ is a penny-farthing we have $y=x_{0}$ because there is a loop at $y$ and $\mu\neq \alpha_{n}$ because there is no loop at $x_{n-1}$. This contradicts lemma \ref{LPF}. Finally, $y$ does not belong to a diamond because of the loop in $y$. Thus theorem 2 holds if one of the contours is a dumb-bell.

Next we look at two penny-farthings. If $x_{0}$ occurs in both we have obviously $x_{0}=x_{0}'$ and $n=n'=2$ by lemma \ref{LPF}.
Then the separated quiver of the full subcategory supported by the two penny-farthings contains a $\tilde{D}_{5}$-diagram. In the remainig case we have $x_{i}=x_{j}'$ for some $i,j$ different from $0$. We claim that the quiver of this category is given by the  picture in figure 3.5 and so the separated quiver contains an $\tilde{A}_{5}$- diagram. Indeed, the claim is obvious for $P(x_{0},x_{0}')= 0=P(x_{0},x_{0}')$ and it follows in the other case because $n=n'=2$ by lemma \ref{LPF}. So the theorem holds for two penny-farthings.

\vspace{0.2cm}
\setlength{\unitlength}{0.8cm}
\begin{picture}(15,2)

\put(4,1){\circle*{0.2}}
\put(6,1){\circle*{0.2}}
\put(8,1){\circle*{0.2}}
\put(4.1,1.1){\vector(1,0){1.8}}
\put(6.1,1.1){\vector(1,0){1.8}}
\put(5.9,0.9){\vector(-1,0){1.8}}
\put(7.9,0.9){\vector(-1,0){1.8}}
\put(3.5,1){\circle{1}}
\put(8.5,1){\circle{1}}
\put(4,1.3){\vector(0,-1){0.2}}
\put(8,1.3){\vector(0,-1){0.2}}\put(4,1.5){$x_{0}$}\put(7.5,1.5){$x_{0}'$}\put(5.5,1.5){$x_{i}=x_{j}'$}

\put(5,0){figure 3.5}
\end{picture}\vspace{0.5cm}

Now, let $C$ be a diamond. We show that $x$ does not belong to another contour. First assume $x$ belongs to a penny-farthing $C'$. Then we have $x=x_{i}'$ for some $i\neq 0$. Lemma \ref{LD} shows $\alpha_{i+1}'=\delta$ or $\alpha_{i+1}'=\gamma$. In the first case we obtain $\alpha_{i+2}'=\beta$ from  lemma  \ref{LD} contradicting the uniqueness lemma for the penny-farthing. Thus we have $z=x_{i+1}'$ and $y\notin C'$ again by the uniqueness lemma. From $\lambda\alpha\neq 0$ we get $P(y,x_{0}')\neq 0$ and therefore $n=2$   by lemma \ref{LPF}. But $z=x_{0}'$ is excluded by the loop in $x_{0}'$. Now let $x$ belong to another diamond $C'$. For $x=x'$ we get $C=C'$ from lemma \ref{LD}.  For $x=z'$ we have $z=y'$ because by lemma \ref{LD} $\gamma$ is the only arrow starting at x and ending in a point with three endomorphisms. By the dual of \ref{LD} $\lambda$ factors through $\alpha'$ or $\beta'$. In the first case we get $P(y,x)\neq 0$. In the second case we have $\lambda=\xi\beta'$. From $0\neq \lambda\alpha= \beta'\xi\alpha$ we get $P(y',t')\neq 0$. Finally assume $x=y'$. Then $\kappa$ factors through $\beta'$ or $\alpha'$ by the dual of lemma \ref{LD}. In the first case $\kappa$ factors even through $\delta'\beta'$ and one obtains from $\kappa\gamma\neq 0$ the contradiction $P(y',x')\neq 0$. For $\kappa=\alpha'\xi$ the first diagram in figure 3.6 is cleaving unless we have $\xi= id$. But then the second diagram is cleaving.

\vspace{0.5cm}
\setlength{\unitlength}{0.8cm}
\begin{picture}(15,3)

\put(2,1){\vector(1,0){1}}
\put(3,1){\vector(1,0){1}}
\put(4,1){\vector(1,0){1}}
\put(4,2){\vector(0,-1){1}}
\put(3,1){\vector(0,-1){1}}

\put(2.5,1.5){$\gamma$}
\put(3.5,0.5){$\xi$}
\put(4.5,0.5){$\alpha'$}
\put(2.5,0.5){$\alpha$}
\put(4.5,1.5){$\gamma'$}

\put(8,1){\vector(1,0){1}}
\put(9,1){\vector(1,0){1}}
\put(9,1){\vector(0,-1){1}}
\put(9,2){\vector(0,-1){1}}

\put(8.5,1.5){$\gamma$}
\put(9.5,0.5){$\alpha'$}
\put(8.5,0.5){$\alpha$}
\put(9.5,1.5){$\gamma'$}

\put(6.5,-0.5){figure 3.6}
\end{picture}\vspace{0.5cm}

By duality, also the point $y$ in a diamond does not belong to another contour. So part a) of the theorem is shown completely.\newpage

For part b) we prove first that the point $t$ in a diamond $C$ does not occur in another contour $C'$. Namely  $C'$ is another diamond with $t=t'$ or a penny-farthing with $t=x_{i}$ for some $i\neq 0$. In both cases there is an arrow  $\epsilon$ in $C'$ ending in $t$ and an arrow $\eta$ starting in $t$ such that $\eta\epsilon \neq 0$. From $\epsilon=\delta$ we get $\beta=\eta$ from lemma \ref{LD}. If $C'$ is a diamond we infer the contradiction $C=C'$. If $C'$ is a penny-farthing we obtain a contradiction to the uniqueness lemma. Thus $\epsilon $ and $\delta$ are different. Dually, $\beta$ and $\eta$ are also different. The first diagram in figure 3.7 is cleaving. For $\lambda\beta=\eta \xi$ implies $0\neq \lambda\beta\delta=\xi\eta\delta$ whence $\eta=\beta$ by \ref{LD}. Similarly $\delta\kappa=\epsilon\xi$ is impossible.

 \setlength{\unitlength}{0.8cm}
\begin{picture}(15,4)

\put(4,3){\vector(1,-1){1}}
\put(4,3){\vector(-1,-1){1}}
\put(6,3){\vector(1,-1){1}}
\put(6,3){\vector(-1,-1){1}}
\put(3,2){\vector(-1,-1){1}}
\put(1,2){\vector(1,-1){1}}
\put(3,2){\vector(1,-1){1}}
\put(5,2){\vector(-1,-1){1}}

\put(0.5,1.5){$\beta$}
\put(2,1.5){$\delta$}
\put(2.8,1.5){$\gamma$}
\put(2.8,2.5){$\kappa$}
\put(4,1.5){$\lambda$}
\put(4,2.5){$\alpha$}
\put(5,2.5){$\beta$}
\put(7,2.5){$\eta$}

\put(9,2){\vector(-1,0){1}}
\put(9,2){\vector(1,0){1}}
\put(10,2){\vector(1,0){1}}
\put(10,2){\vector(0,-1){1}}
\put(9,1){\vector(1,0){1}}
\put(9,2){\vector(0,-1){1}}
\put(8,1){\vector(1,0){1}}
\put(10.5,1.5){$\gamma$}
\put(9.5,0.5){$\lambda$}
\put(8.5,0.5){$\beta$}
\put(9.5,2.2){$\kappa$}
\put(10.5,2.2){$\delta$}

\put(8.5,2.2){$\alpha '$}

\put(6.5,-0.5){figure 3.7}\end{picture}\vspace{0.5cm}

We now look at the case where a diamond and a penny-farthing intersect at $z=x_{i}'$ for some $i \neq 0$. We claim that up to duality the full  subcategory supported by $x,y,z,x_{0}'$ has one of the forms given in figure  3.8. One verifies that in both cases the separated quiver contains a $\tilde{D}_{n}$-diagram. Clearly, the first form occurs if $P(x_{0}',a)=0=P(a,x_{0}')$ holds for all $a$ outside of the penny-farthing. If this is not true we have $n=2$ and $\gamma, \alpha, \rho', \alpha_{1}', \alpha_{2}'$ all remain arrows. Suppose we have $\lambda= \xi\lambda'$ for some arrow $\lambda'$ ending in $x_{0}$. Then $\kappa$ remains an arrow, because otherwise there are three arrows starting and three arrows ending in $x_{0}'$.

\vspace{0.5cm}
\setlength{\unitlength}{0.8cm}
\begin{picture}(15,6)

\put(1,1){\circle*{0.2}}\put(0.5,1){$x$}
\put(3,3){\circle*{0.2}}\put(3.5,4.8){$x_{0}'$}
\put(3,5){\circle*{0.2}}\put(3.5,3){$z$}
\put(5,1){\circle*{0.2}}\put(5.2,1){$y$}
\put(1,1.2){\vector(1,1){1.8}}
\put(3,2.8){\vector(-1,-1){1.8}}
\put(2.9,3.1){\vector(0,1){1.8}}
\put(3.1,4.9){\vector(0,-1){1.8}}
\put(3,5.5){\circle{1}}
\put(2.9,3.1){\vector(0,1){1.8}}
\put(4.8,1){\vector(-1,1){1.8}}
\put(3.2,2.9){\vector(1,-1){1.8}}

\put(3.25,5){\vector(-1,0){0.2}}

\put(8,1){\circle*{0.2}}\put(7.5,1){$x$}
\put(10,3){\circle*{0.2}}\put(9.3,3){$z$}
\put(10,5){\circle*{0.2}}\put(9.2,4.8){$x_{0}'$}
\put(12,1){\circle*{0.2}}\put(12.2,1){$y$}
\put(8,1.2){\vector(1,1){1.8}}
\put(10,2.8){\vector(-1,-1){1.8}}
\put(9.9,3.1){\vector(0,1){1.8}}
\put(10.1,4.9){\vector(0,-1){1.8}}
\put(10,5.5){\circle{1}}
\put(9.9,3.1){\vector(0,1){1.8}}
\put(12,1.2){\vector(-1,2){1.8}}
\put(10.2,2.9){\vector(1,-1){1.8}}

\put(10.25,5){\vector(-1,0){0.2}}

\put(5.5,0){figure 3.8}
\end{picture}\vspace{0.3cm}

The only case still not treated is when two diamonds meet only in $z=z'$. Then we consider the second diagram in figure 3.7. If this is not cleaving we have $\gamma\kappa=\xi\alpha'$ or $\delta\kappa=\xi\alpha'$, whence we get $\alpha'\gamma\neq 0$ after multiplication with $\gamma$. Using symmetry one also gets $\alpha\gamma' \neq 0$. One finds an obvious cleaving diagram of type $\tilde{D}_{4}$ with arrows $\gamma,\gamma',\alpha, \alpha'$.

\section{Mild contours and minimal representation-infinite ray categories}

\subsection{Dumb-bells}

\begin{proposition}\label{MIDB} Let $C$ be a dumb-bell in a ray category $P$ such that $P/\rho^{2}\mu$ is representation-finite. Assume that $\tau:y \rightarrow z$ is an additional arrow. Then we have:
\begin{enumerate}
 \item $\tau$ and $\rho$ are the only arrows starting in $y$.
\item $\lambda$ and $\mu$ are the only arrows starting or ending in $x$.
\item $\lambda^{3}=0$ and $\rho^{3}=0$.
\item If $\tau\mu \neq 0$ there is no arrow ending in $z$ except $\tau$ and $\zeta\tau=0$ holds for all arrows starting in $z$.
\end{enumerate}

\end{proposition}

 Proof: a) If there is another arrow $\tau'$ starting in $y$ the first diagram in figure 4.1 is cleaving. Namely for $\rho^{2}=\tau\xi$ the morphism  $\xi$ is of the form $\rho\xi'$ or $\mu\xi'$ by \ref{LDB}. The first case implies $\rho=\xi'\tau$ and the second gives the non-zero morphism $\xi'\tau$ in $P(y,x)$.

b) An additional arrow $\zeta:t \rightarrow x$ gives rise to the second diagram in figure 4.1 which is cleaving by the arguments used in a) and their duals.

\vspace{0.5cm}
\setlength{\unitlength}{0.8cm}
\begin{picture}(15,5)

\put(1,2){\vector(1,0){1}}\put(2,3){\vector(1,0){1}}\put(4,1){\vector(1,0){1}}\put(5,2){\vector(1,0){1}}\put(5,3){\vector(1,0){1}}\put(6,4){\vector(1,0){1}}\put(9,3){\vector(1,0){1}}\put(9,4){\vector(1,0){1}}\put(10,5){\vector(1,0){1}}\put(12,2){\vector(1,0){1}}\put(13,3){\vector(1,0){1}}

\put(2,2){\vector(0,-1){1}}\put(2,3){\vector(0,-1){1}}\put(5,2){\vector(0,-1){1}}\put(5,3){\vector(0,-1){1}}\put(6,3){\vector(0,-1){1}}\put(6,4){\vector(0,-1){1}}\put(9,2){\vector(0,-1){1}}\put(9,3){\vector(0,-1){1}}\put(9,4){\vector(0,-1){1}}\put(10,4){\vector(0,-1){1}}\put(10,5){\vector(0,-1){1}}\put(13,2){\vector(0,-1){1}}\put(13,3){\vector(0,-1){1}}\put(13,4){\vector(0,-1){1}}

\put(2,3){\vector(0,1){1}}

\put(2.2,1.5){$\rho$}\put(2.2,2.5){$\rho$}\put(6.2,2.5){$\rho$}\put(6.2,3.5){$\rho$}\put(10.2,3.5){$\rho$}\put(10.2,4.5){$\rho$}\put(13.2,1.5){$\rho$}\put(13.2,2.5){$\rho$}\put(13.2,3.5){$\rho$}
\put(4.5,1.5){$\lambda$}\put(4.5,2.5){$\lambda$}\put(8.5,1.5){$\lambda$}\put(8.5,2.5){$\lambda$}\put(8.5,3.5){$\lambda$}
\put(1.5,2.2){$\mu$}\put(5.5,1.5){$\mu$}\put(5.5,3.2){$\mu$}\put(9.5,2.5){$\mu$}\put(9.5,4.2){$\mu$}\put(12.5,2.2){$\mu$}
\put(2.5,3.2){$\tau$}\put(6.5,4.2){$\tau$}\put(10.5,5.2){$\tau$}\put(13.5,3.2){$\tau$}\put(4.5,0.5){$\zeta$}\put(1.5,3.5){$\tau'$}
\put(5.5,0){figure 4.1}
\end{picture}\vspace{0.5cm}
c) $\lambda^{3}=0$ follows from the third diagram in figure 4.1 which is cleaving again. So assume $\rho^{3}\neq 0$. Then the fourth diagram in figure 4.1 is cleaving for $\tau\rho \neq 0$ and the first in figure 4.2 for  $\tau\rho = 0$.

d) For $\tau\mu \neq 0$ an obvious $\tilde{D}_{4}$- diagram shows $\tau\rho=0.$ Then the second and the third diagram in figure 4.2 show that there is no additional arrow $\zeta$ ending in $z$ and that the composition of $\tau$ with all arrows $\zeta$ is zero.

\vspace{0.5cm}
\setlength{\unitlength}{0.8cm}
\begin{picture}(15,5)

\put(1,2){\vector(0,-1){1}}\put(1,3){\vector(0,-1){1}}\put(2,3){\vector(0,-1){1}}\put(2,4){\vector(0,-1){1}}\put(2,5){\vector(0,-1){1}}
\put(1,2){\vector(1,0){1}}\put(1,3){\vector(1,0){1}}\put(2,4){\vector(1,0){1}}\put(2,5){\vector(1,0){1}}

\put(5,2){\vector(0,-1){1}}\put(5,3){\vector(0,-1){1}}\put(6,3){\vector(0,-1){1}}\put(6,4){\vector(0,-1){1}}
\put(5,2){\vector(1,0){1}}\put(5,3){\vector(1,0){1}}\put(6,3){\vector(1,0){1}}\put(6,4){\vector(1,0){1}}\put(8,3){\vector(-1,0){1}}

\put(10,2){\vector(0,-1){1}}\put(10,3){\vector(0,-1){1}}\put(11,3){\vector(0,-1){1}}\put(11,4){\vector(0,-1){1}}
\put(10,2){\vector(1,0){1}}\put(10,3){\vector(1,0){1}}\put(11,3){\vector(1,0){1}}\put(11,4){\vector(1,0){1}}\put(12,4){\vector(1,0){1}}

\put(0.5,1.5){$\lambda$}\put(0.5,2.5){$\lambda$}\put(4.5,1.5){$\lambda$}\put(4.5,2.5){$\lambda$}\put(9.5,1.5){$\lambda$}\put(9.5,2.5){$\lambda$}
\put(2.2,2.5){$\rho$}\put(2.2,3.5){$\rho$}\put(2.2,4.5){$\rho$}\put(6.2,2.5){$\rho$}\put(6.2,3.5){$\rho$}\put(11.2,2.5){$\rho$}\put(11.2,3.5){$\rho$}
\put(1.5,1.5){$\mu$}\put(1.5,3.2){$\mu$}\put(5.5,1.5){$\mu$}\put(5.5,3.2){$\mu$}\put(10.5,1.5){$\mu$}\put(10.5,3.2){$\mu$}

\put(2.7,4.2){$\tau$}\put(2.7,5.2){$\tau$}\put(6.7,4.2){$\tau$}\put(6.7,3.2){$\tau$}\put(11.7,4.2){$\tau$}\put(11.7,3.2){$\tau$}
\put(7.5,3.2){$\zeta$}\put(12.5,4.2){$\zeta$}
\put(5,0.5){$\tau\rho=0$ in all three cases}
\put(5,-0.5){figure 4.2}
\end{picture}\vspace{0.5cm}
\subsection{Penny farthings} For the convenience of the reader we give a full proof of the following result from \cite{Treue} which was proven there using coverings and a lemma  from \cite{BGRS}. The proof is not as easy as claimed in \cite[section 13.15]{BGRS}.
\begin{proposition}\label{MIPF} Let $C$ be a penny-farthing in a ray category $P$ such that $P'/\rho^{3}$ is representation-finite for each full subcategory $P'$ which contains $x_{0}$ and at most four other points. Suppose that $P(x_{0},y)\neq0 $ for some $y$ not in $C$. Then we have $n=2$ and we are in one of the following three situations:
\begin{enumerate}
 \item There is an arrow $\beta:x_{0} \rightarrow b$ 
and this is the only arrow outside $C$ ending or starting in $C$. We have $\beta\rho=0$,$\delta\beta=0$ for all arrows $\delta$ starting in $b$ and $y=b$.
\item There is an arrow $\gamma:x_{1} \rightarrow c$ and this is the only arrow outside $C$ ending or starting in $C$. We have $\gamma\alpha_{1}\rho=0$,$\delta\gamma=0$ for all arrows starting in  $c$ and $y=c$.
\item There is an arrow $\beta:x_{0} \rightarrow b$ as well as an arrow $\gamma:x_{1}  \rightarrow c$. These are the only two arrows outside $C$ ending or starting in $C$. We have $0=\beta\rho=\beta\alpha_{2}$,$0=\delta\beta$ for all arrows starting in $b$,  $0=\gamma\alpha_{1}$,$0=\epsilon\gamma$ for all arrows starting in $c$ and $y=b\neq c$.
\end{enumerate}
In all three cases there are no additional arrows ending in $b$ or $c$.

\end{proposition}

Proof: We know already that $n=2$ from lemma \ref{LPF}. The relations $\alpha_{1}\rho=\eta\alpha_{1}$ and $\rho\alpha_{2}=\alpha_{2}\eta$ are excluded because both  imply $\eta=\alpha_{1}\rho\alpha_{2}$  and the final contradictions $0=\alpha_{1}\rho$ and $0=\rho\alpha_{2}$. 

First we treat the case that there is an arrow $\beta:x_{0} \rightarrow b$. Then the separated quiver shows that there is no additional arrow ending at $x_{0}$ and we get $\beta\rho=0$ from an obvious $\tilde{D}_{4}$-diagram. Furthermore $\rho^{2}=\eta\beta$ would imply $0\neq \rho^{3}=\eta\beta\rho=0$.

If $\delta:z \rightarrow x_{1}$ is another arrow the separated quiver shows that only $\delta$ and $\alpha_{1}$ end at $x_{1}$. In the case $\alpha_{2}\delta \neq 0$ the first diagram in figure 4.3 is cleaving. Namely, from $\alpha_{2}\delta=\rho\eta$ we get $\eta=\rho\eta'$ whence the contradiction $0\neq \alpha_{2}\delta=\rho^{2}\eta'=\alpha_{2}\alpha_{1}\eta'$ or $\eta=\alpha_{2}\eta'$ with $\eta'=\delta\eta''$ or $\eta'=\eta''\delta$ which are both impossible. For $\alpha_{2}\delta =0$ the second diagram in figure 4.3 is cleaving. Namely $\alpha_{1}\rho=\delta\eta$ implies $0\neq \rho^{3}=\alpha_{2}\alpha_{1}\rho=\alpha_{2}\delta\eta=0$.  From $\alpha_{1}\rho=\eta\beta$ one gets $\eta=\alpha_{1}\eta'$ or $\eta=\delta\eta'$ which imply $\rho=\eta'\beta$ or $0\neq \alpha_{2}\alpha_{1}\rho=\alpha_{2}\delta\eta'\beta=0$.

\begin{picture}(15,4)

\put(1,3){\vector(1,0){1}}\put(2,3){\vector(0,-1){1}}\put(3,3){\vector(0,-1){1}}\put(3,3){\vector(-1,0){1}}\put(3,3){\vector(1,0){1}}
\put(3,2){\vector(-1,0){1}}\put(4,2){\vector(-1,0){1}}

\put(7,3){\vector(1,0){1}}\put(8,3){\vector(0,-1){1}}\put(9,3){\vector(0,-1){1}}\put(9,3){\vector(1,0){1}}
\put(9,3){\vector(-1,0){1}}\put(9,2){\vector(-1,0){1}}\put(9,2){\vector(0,-1){1}}\put(10,2){\vector(-1,0){1}}\put(10,1){\vector(-1,0){1}}

\put(3.2,2.5){$\rho$}\put(8.7,2.5){$\rho$}\put(2.5,1.5){$\rho$}\put(8.5,1.7){$\rho$}
\put(2.5,3.2){$\alpha_{1}$}\put(8.5,3.2){$\alpha_{1}$}\put(9.2,1.5){$\alpha_{1}$}
\put(1.5,2.5){$\alpha_{2}$}\put(3.5,1.5){$\alpha_{2}$}\put(7.5,2.5){$\alpha_{2}$}\put(9.5,2.2){$\alpha_{2}$}
\put(3.5,3.2){$\beta$}\put(9.5,3.2){$\beta$}
\put(1.5,3.2){$\delta$}\put(7.5,3.2){$\delta$}\put(9.5,0.5){$\delta$}
\put(11,2){$\alpha_{2}\delta=0=\alpha_{1}\alpha_{2}$}

\put(5,-0.5){figure 4.3}
\end{picture}\vspace{0.5cm}

Thus we have shown that there is no additional arrow ending in $x_{1}$. This makes the relation $\eta\beta=\alpha_{1}\rho\neq 0$ impossible because of $\eta=\alpha_{1}\eta'$ whence $\rho=\eta'\beta$. 

For an arrow $\epsilon:z \rightarrow b$ different from $\beta$ the first diagram in figure 4.4 is obviously cleaving and for an arrow $\epsilon:b \rightarrow z$  with $\epsilon\beta \neq 0$ the second diagram. Namely, only factorizations $0\neq\epsilon\beta=\zeta\rho$ and $\epsilon\beta=\eta\alpha_{1}$ have to be excluded. Writing in the first case $\zeta=\zeta'\rho^{i}$ for some $i\geq 1$ such that $\zeta'$ does not factor through $\rho$ we get $\zeta=\eta\alpha_{1}\rho$ from $\beta\rho=0=\alpha_{1}\rho^{2}$ for some $\eta$ ( which can be different from the one considered in the second factorization ).

We claim that in both cases the third diagram in figure 4.4 is cleaving.
For $\eta=\epsilon\xi$ implies the contradiction $0\neq \epsilon\beta=\epsilon\xi\alpha_{1}\rho$ in the first case and $0\neq \epsilon\beta=\epsilon\xi\alpha_{1}$ in the second. Similarly $\eta=\xi\alpha_{2}$ gives $0\neq \epsilon\beta=\xi\alpha_{2}\alpha_{1}\rho=\xi\rho^{3}$ or $0\neq \epsilon\beta=\xi\alpha_{2}\alpha_{1}=\xi\rho^{2}$. Now $\phi\rho^{3}=0$ holds for any arrow $\phi$ and $\alpha_{1}\rho^{2}=0$. Thus in the second case $\xi$ starts with $\rho$ and the right hand sides are $0$ in both cases. Finally we look at $\rho\alpha_{2}=\xi\eta$. Then $\xi$  cannot end with $\rho$ whence it ends with $\alpha_{2}$. Since $\alpha_{1}$is the only arrow ending at $x_{1}$ we obtain  $0\neq \rho\alpha_{2}=\alpha_{2}\alpha_{1}\xi'\eta=\rho^{2}\xi'\eta$ contradicting the irreducibility of $\alpha_{2}$.

\begin{picture}(20,4)

\put(0,2){\vector(1,0){1}}
\put(2,2){\vector(0,-1){1}}
\put(2,2){\vector(1,0){1}}
\put(3,2){\vector(1,0){1}}
\put(2,2){\vector(-1,0){1}}
\put(3,3){\vector(0,-1){1}}
\put(2,1){\vector(1,0){1}}
\put(3,2){\vector(0,-1){1}}

\put(6,2){\vector(-1,0){1}}
\put(7,2){\vector(0,-1){1}}
\put(7,2){\vector(1,0){1}}
\put(8,2){\vector(1,0){1}}
\put(7,2){\vector(-1,0){1}}
\put(8,3){\vector(0,-1){1}}
\put(7,1){\vector(1,0){1}}
\put(8,2){\vector(0,-1){1}}

\put(11,2){\vector(-1,0){1}}
\put(11,2){\vector(0,-1){1}}
\put(11,2){\vector(1,0){1}}
\put(12,2){\vector(0,-1){1}}
\put(13,2){\vector(-1,0){1}}
\put(13,2){\vector(1,0){1}}
\put(15,2){\vector(-1,0){1}}
\put(12,2){\vector(0,1){1}}
\put(11,1){\vector(1,0){1}}

\put(2.5,2.2){$\rho$}
\put(3.2,1.5){$\rho$}
\put(7.5,2.2){$\rho$}
\put(8.2,1.5){$\rho$}
\put(11.5,2.2){$\rho$}
\put(12.2,1.5){$\rho$}
\put(1.5,1.5){$\alpha_{1}$}
\put(3.5,2.2){$\alpha_{1}$}
\put(6.5,1.5){$\alpha_{1}$}
\put(8.5,2.2){$\alpha_{1}$}
\put(10.5,1.5){$\alpha_{1}$}
\put(12.2,2.5){$\alpha_{1}$}

\put(1.5,2.2){$\beta$}\put(6.5,2.2){$\beta$}\put(10.5,2.2){$\beta$}
\put(0.5,2.2){$\epsilon$}
\put(14.5,2.2){$\epsilon$}
\put(5.5,2.2){$\epsilon$}

\put(2.5,0.5){$\alpha_{2}$}
\put(3.2,2.7){$\alpha_{2}$}
\put(8.2,2.7){$\alpha_{2}$}
\put(7.5,0.5){$\alpha_{2}$}
\put(11.5,0.5){$\alpha_{2}$}
\put(12.5,2.2){$\alpha_{2}$}
\put(13.5,2.2){$\eta$}
\put(13,2){\circle*{0.1}}

\put(4,0){$0=\alpha_{1}\alpha_{2}$ in all cases}

\put(5,-1){figure 4.4}
\end{picture}\vspace{1cm}

Assume now that there is an arrow $\gamma:x_{1} \rightarrow c$ different from $\alpha_{2}$. Then there is no third arrow $\gamma'$ starting in $x_{1}$. For otherwise the separated quiver shows that only $\rho$ and $\alpha_{1}$ start at $x_{0}$ and then the first diagram in figure 4.5 is cleaving because $0\neq \rho\alpha_{2}= \eta\gamma$ equals $\rho\eta'\gamma$ or $\alpha_{2}\eta'\gamma=\alpha_{2}\alpha_{1}\rho\alpha_{2}=0$. 

For $\gamma\alpha_{1}\neq 0$ the second diagram is cleaving. Namely the relation $\gamma\alpha_{1}=\eta\beta$ is impossible because we already know $\delta\beta=0$ for all arrows $\delta$. From $\gamma\alpha_{1}=\eta\rho$ we get $\gamma\alpha_{1}=\gamma\alpha_{1}\rho$ since $0=\beta\rho=\alpha_{1}\rho^{2}$ and only $\rho$,$\beta$,$\alpha_{1}$ start at $x_{0}$ and since only $\gamma$,$\alpha_{2}$ start at $x_{1}$ and $\phi\rho^{3}=0$ for all arrows.

Next, in case $\beta\alpha_{2}\neq 0$ the third diagram of figure 4.5 is cleaving because  $\rho\alpha_{2}=\eta\gamma$ gives $0\neq \rho^{3}=\eta\gamma\alpha_{1}=0$.

Finally, if $\epsilon\gamma\neq 0$ holds for some $\epsilon$ we consider the last diagram in figure 4.5. From $\epsilon\gamma=\eta \alpha_{2}$ we obtain $\epsilon\gamma=\xi\epsilon\gamma\alpha_{1}\rho\alpha_{2}$ or $\epsilon\gamma=\epsilon\xi\gamma\alpha_{1}\rho\alpha_{2}$ which are both impossible.

\begin{picture}(20,3.5)

\put(1,1){\vector(-1,0){1}}
\put(1,1){\vector(0,1){1}}
\put(1,1){\vector(1,0){1}}
\put(2,2){\vector(0,-1){1}}
\put(2,1){\vector(1,0){1}}

\put(6,1){\vector(-1,0){1}}
\put(6,1){\vector(1,0){1}}
\put(6,2){\vector(1,0){1}}
\put(6,2){\vector(-1,0){1}}
\put(6,2){\vector(0,-1){1}}
\put(7,2){\vector(0,-1){1}}
\put(8,2){\vector(-1,0){1}}

\put(11,2){\vector(-1,0){1}}
\put(11,2){\vector(0,-1){1}}
\put(11,2){\vector(1,0){1}}
\put(11,1){\vector(1,0){1}}
\put(12,2){\vector(-1,1){1}}
\put(12,2){\vector(1,0){1}}
\put(12,2){\vector(0,-1){1}}
\put(12,3){\vector(0,-1){1}}
\put(12,3){\vector(1,0){1}}

\put(6,-2){\vector(-1,0){1}}\put(6.5,-1.8){$\rho$}\put(5.5,-1.8){$\beta$}
\put(6.5,-3.4){$\alpha_{2}$}\put(7.2,-1.5){$\alpha_{2}$}\put(7.2,-2.5){$\rho$}
\put(5.5,-2.5){$\alpha_{1}$}\put(7.8,-2.4){$\alpha_{1}$}\put(7.5,-0.8){$\gamma$}\put(8.5,-0.8){$\epsilon$}
\put(9,-2){$0=\alpha_{1}\alpha_{2}$}
\put(6,-2){\vector(0,-1){1}}
\put(6,-2){\vector(1,0){1}}
\put(6,-3){\vector(1,0){1}}
\put(8,-1){\vector(1,0){1}}
\put(7,-2){\vector(1,0){1}}
\put(7,-2){\vector(0,-1){1}}
\put(7,-1){\vector(0,-1){1}}
\put(7,-1){\vector(1,0){1}}

\put(2.2,1.5){$\rho$}
\put(2.5,0.5){$\rho$}
\put(6.5,2.2){$\rho$}
\put(7.2,1.5){$\rho$}
\put(11.5,1.7){$\rho$}
\put(12.2,1.5){$\rho$}
\put(5.5,1.5){$\alpha_{1}$}
\put(10.5,1.5){$\alpha_{1}$}
\put(12.2,2.5){$\alpha_{2}$}

\put(5.5,2.2){$\beta$}
\put(11.5,3){$\beta$}
\put(10.5,2.2){$\beta$}

\put(1.5,0.5){$\alpha_{2}$}
\put(6.5,0.5){$\alpha_{2}$}
\put(7.5,2.2){$\alpha_{2}$}
\put(11.5,0.5){$\alpha_{2}$}
\put(12.5,2.2){$\alpha_{1}$}
\put(0.5,0.5){$\gamma$}
\put(5.5,0.5){$\gamma$}
\put(12.5,3.2){$\gamma$}
\put(0.5,1.5){$\gamma'$}

\put(13,1.5){$0=\alpha_{1}\alpha_{2}$}\put(13,1){$0=\beta\rho$}

\put(5,-4){figure 4.5}
\end{picture}\vspace{4cm}

Using all the information obtained so far the reader can easily check that we are in case a) or c) of the lemma if there is an arrow $\beta:x_{0}\rightarrow b$.

So assume there is no such arrow. Then the  assumption $P(x_{0},y)\neq 0$ implies the existence of an arrow $\gamma:x_{1} \rightarrow c$ with $\gamma\alpha_{1}\neq 0$ and $\gamma$ and $\alpha_{2}$ are the only arrows starting at $x_{1}$ as shown before. The relation $\gamma\alpha_{1}=\eta\rho=\eta'\alpha_{1}\rho$ implies $\eta'=\xi\gamma$ or $\eta'=\gamma\xi$ which are both impossible.

If $\delta$ is another arrow ending in $x_{0}$, the first diagram in figure 4.6 is cleaving provided $\rho^{2}=\delta\eta$ is excluded. But such a factorization contradicts the uniqueness lemma \ref{wegeindeutig}. Similarly for another arrow $\delta$ ending at $x_{1}$ the second diagram in figure 4.6 is cleaving because in a relation $\alpha_{1}\rho=\delta\eta$ the morphism $\eta$ starts with $\rho$ or with $\alpha_{1}$ leading always to a contradiction. If we  have $\delta\gamma\neq 0$ for some  $\delta$ the third diagram is cleaving because the relation $\epsilon\gamma=\eta\alpha_{2}=\eta'\gamma\alpha_{1}\rho\alpha_{2}$ with $\eta'=\epsilon\xi$ or $\eta'=\xi\epsilon$ is excluded. Finally for $\gamma\alpha_{1}\rho\neq 0$ we look at the last diagram. A relation $\gamma\alpha_{1}\rho=\eta\alpha_{1}$ implies $\eta\alpha_{1}=\eta'\gamma\alpha_{1}$ and this term vanishes by the fact shown just before.  The proof of the proposition is complete.

\vspace{1cm}
\begin{picture}(15,8)

\put(2,4){\vector(-1,0){1}}\put(2.5,5.2){$\rho$}\put(3.2,4.6){$\rho$}
\put(2,4){\vector(1,0){1}}
\put(4,4){\vector(-1,0){1}}\put(1.5,4.6){$\alpha_{1}$}\put(3.5,5.2){$\alpha_{2}$}\put(2.3,4.1){$\alpha_{2}$}
\put(2,5){\vector(0,-1){1}}\put(1.4,4.1){$\gamma$}\put(3.5,4.1){$\delta$}
\put(2,5){\vector(1,0){1}}
\put(3,5){\vector(0,-1){1}}
\put(4,5){\vector(-1,0){1}}

\put(8,4){\vector(-1,0){1}}\put(9.2,5.5){$\alpha_{2}$}\put(8.3,4.1){$\alpha_{2}$}
\put(8,4){\vector(1,0){1}}\put(2.5,5.2){$\rho$}\put(3.2,4.6){$\rho$}
\put(8,5){\vector(1,0){1}}\put(7.5,4.6){$\alpha_{1}$}\put(9.2,5.1){$\alpha_{1}$}
\put(9,5){\vector(1,0){1}}\put(10.5,5.1){$\delta$}
\put(9,5){\vector(0,-1){1}}\put(8.5,5.2){$\rho$}\put(9.2,4.6){$\rho$}
\put(11,5){\vector(-1,0){1}}\put(7.4,4.1){$\gamma$}\put(9.4,6.1){$\gamma$}
\put(9,6){\vector(1,0){1}}
\put(9,6){\vector(0,-1){1}}
\put(8,5){\vector(0,-1){1}}

\put(8,1){\vector(-1,0){1}}\put(8.5,2.2){$\rho$}\put(9.2,1.6){$\rho$}
\put(8,1){\vector(1,0){1}}\put(0.4,1.1){$\gamma$}\put(2.4,3.1){$\gamma$}
\put(8,2){\vector(1,0){1}}\put(7.5,1.6){$\alpha_{1}$}\put(9.2,2.1){$\alpha_{1}$}
\put(9,2){\vector(1,0){1}}\put(9.2,2.5){$\alpha_{2}$}\put(8.3,1.1){$\alpha_{2}$}
\put(9,2){\vector(0,-1){1}}
\put(10,2){\vector(1,0){1}}
\put(9,3){\vector(1,0){1}}
\put(9,3){\vector(0,-1){1}}
\put(8,2){\vector(0,-1){1}}

\put(1,1){\vector(-1,0){1}}\put(1.5,2.2){$\rho$}\put(2.2,1.6){$\rho$}
\put(1,1){\vector(1,0){1}}\put(7.4,1.1){$\gamma$}\put(9.4,3.1){$\gamma$}
\put(1,2){\vector(1,0){1}}\put(3.4,3.1){$\delta$}
\put(2,2){\vector(1,0){1}}\put(0.5,1.6){$\alpha_{1}$}\put(2.2,2.1){$\alpha_{1}$}
\put(2,2){\vector(0,-1){1}}\put(2.2,2.5){$\alpha_{2}$}\put(1.3,1.1){$\alpha_{2}$}
\put(3,3){\vector(1,0){1}}\put(10.4,2.1){$\gamma$}
\put(2,3){\vector(1,0){1}}\put(3.5,0){$\alpha_{1}\alpha_{2}=0$ in the last three diagrams}
\put(2,3){\vector(0,-1){1}}
\put(1,2){\vector(0,-1){1}}

\put(5,-1){figure 4.6}
\end{picture}\vspace{1cm}

\subsection{Diamonds}
\begin{proposition}\label{MID}
 Let $D=\{\alpha\gamma,\beta\delta\}$ be a diamond. Write $\lambda$ and $\kappa$ as products of irreducible morphisms as in figure 4.7. Then we have:
\begin{enumerate}
 \item The diagram of figure 4.7 is cleaving in $P$.
\item The decompositions of $\lambda$ and $\kappa$ are unique and $s$ and $t$ are at most 2.
\item The figure contains all the non-zero paths starting at $x$ except for $s=1$ and $r=2$ where there can be the additional morphism $\kappa_{1}\beta\kappa\alpha\gamma$.
\item The projective indecomposable $k(P)$-module $k(P)(x,-)$ is also injective.
\end{enumerate}

\end{proposition}

\vspace{0.5cm}
\setlength{\unitlength}{0.8cm}
\begin{picture}(10,4)

\put(1,2){\vector(1,1){1}}
\put(1,2){\vector(1,-1){1}}
\put(2,1){\vector(1,1){1}}
\put(2,3){\vector(1,-1){1}}
\put(2,3){\vector(2,0){2}}
\put(3,2){\vector(2,0){2}}
\put(7,3){\vector(2,0){2}}
\put(8,2){\vector(2,0){2}}

\put(9,3){\vector(1,1){1}}
\put(9,3){\vector(1,-1){1}}
\multiput(4.1,3)(0.2,0){15}{\circle*{0.05}}
\multiput(5.1,2)(0.2,0){15}{\circle*{0.05}}
\put(10,2){\vector(1,1){1}}
\put(10,4){\vector(1,-1){1}}
\put(5,0){figure 4.7}

\put(1.2,2.5){$\gamma$}
\put(2.8,2.5){$\alpha$}
\put(8.8,2.3){$\gamma$}
\put(10.8,2.3){$\alpha$}
\put(1.2,1.2){$\delta$}
\put(2.8,1.2){$\beta$}\put(3.8,1.5){$\kappa_{1}$}\put(9,1.5){$\kappa_{t}$}
\put(8.8,3.4){$\delta$}
\put(2.8,3.4){$\lambda_{1}$}
\put(7.8,3.4){$\lambda_{s}$}
\put(10.8,3.4){$\beta$}

\end{picture}\vspace{0.5cm}

Proof: a) For $t=0$ only the possibility $\xi\kappa=\gamma\lambda_{s}\ldots\lambda_{2}$ has to be excluded. But this would imply the contradiction
$\alpha=\xi\lambda_{1}$. Thus the diagram is cleaving for $t=0$. As already observed this implies $s\leq 2$ and dually $t\leq 2$. For $s=t=2$, up to duality only the case $\kappa_{1}\alpha=\xi\lambda_{1}$ has to be considered. But then $P'$ contains a category with number 85 form the BHV-list as a cleaving diagram as shown as the first picture in figure 4.8. Here again up to duality only $\kappa_{1}\beta=\xi\eta$ has to be excluded. But this induces the cleaving diagram $\tilde{A}(\beta,\eta,\lambda_{1},\alpha)$ as one easily verifies. 

b) This is clear from lemma \ref{wegeindeutig} and the discussion above.

c) We know already from lemma \ref{LD} that another arrow $\phi$  starting at $x$  or satisfying $\phi\delta\neq 0$ does not exist.

\vspace{0.5cm}
\setlength{\unitlength}{0.8cm}
\begin{picture}(20,4)

\put(0,2){\vector(1,0){1}}\put(1.4,3.2){$\lambda_{1}$}\put(2.4,3.2){$\lambda_{2}$}\put(3.4,3.2){$\delta$}
\put(1,2){\vector(1,0){1}}\put(1.4,1.5){$\kappa_{1}$}\put(2.4,1.5){$\kappa_{2}$}\put(0.5,1.5){$\beta$}
\put(2,2){\vector(1,0){1}}\put(0.6,2.5){$\alpha$}\put(1.6,2.5){$\xi$}\put(2.6,2.5){$\gamma$}
\put(1,3){\vector(1,0){1}}
\put(2,3){\vector(1,0){1}}
\put(3,3){\vector(1,0){1}}
\put(1,3){\vector(0,-1){1}}
\put(2,3){\vector(0,-1){1}}
\put(3,3){\vector(0,-1){1}}

\put(8,3){\vector(1,0){1}}\put(8.4,3.2){$\lambda_{2}$}\put(8.4,1.5){$\kappa_{2}$}\put(9.5,1.5){$\alpha$}\put(10.4,1.5){$\kappa_{1}$}\put(9.5,3.2){$\delta$}\put(8.6,2.5){$\gamma$}\put(9.6,2.5){$\beta$}
\put(8,2){\vector(1,0){1}}
\put(9,2){\vector(1,0){1}}
\put(9,3){\vector(1,0){1}}
\put(9,3){\vector(0,-1){1}}
\put(10,3){\vector(0,-1){1}}

\put(10,2){\vector(1,0){1}}

\put(5,0){figure 4.8}

\end{picture}\vspace{0.5cm}

The obvious $\tilde{D}_{4}$-diagrams make an arrow $\phi$ with $\phi\gamma \neq 0$ or with $\phi\beta\delta \neq 0$ or with $0 \neq \phi\kappa_{1}\beta\delta$ impossible. For $0 \neq \phi\lambda_{1}\gamma$ the first diagram in figure 4.9 is cleaving.
In the same vein an additional arrow $\phi$ with $\phi\kappa\beta\delta\neq 0$ or one with $\phi\lambda\gamma\neq 0$ gives obvious $\tilde{D}_{4}$-diagrams. 

Finally for $\phi$ with $\phi\beta\delta\lambda\gamma \neq 0$
another $\tilde{D}_{4}$-diagram implies $\phi=\kappa_{1}$. Because the second diagram of figure 4.8 cannot be cleaving in $P'$ we have $s=1$. The last non-zero path cannot be prolonged, because  $\kappa_{2}$ is the only candidate for an arrow and $\kappa\alpha\gamma\lambda\gamma=0$.

\vspace{0.5cm}
\setlength{\unitlength}{0.8cm}
\begin{picture}(20,4)

\put(7,2){\vector(1,0){1}}
\put(7.4,1.5){$\kappa_{1}$}
\put(8.4,1.5){$\kappa_{2}$}
\put(9.5,1.5){$\alpha$}
\put(10.4,1.5){$\kappa_{1}$}
\put(11.4,1.5){$\phi$}
\put(9.5,3.2){$\delta$}
\put(8.6,2.5){$\gamma$}
\put(9.6,2.5){$\beta$}
\put(8,2){\vector(1,0){1}}
\put(9,2){\vector(1,0){1}}
\put(9,3){\vector(1,0){1}}
\put(9,3){\vector(0,-1){1}}
\put(10,3){\vector(0,-1){1}}

\put(10,2){\vector(1,0){1}}
\put(12,2){\vector(-1,0){1}}

\put(0,2){\vector(1,0){1}}
\put(2.2,2.2){$\phi$}
\put(1.4,2.2){$\lambda_{1}$}\put(0.5,1.5){$\alpha$}
\put(0.5,2.2){$\gamma$}
\put(2.4,1.5){$\lambda_{2}$}

\put(1,2){\vector(1,0){1}}
\put(2,2){\vector(1,0){1}}
\put(1,2){\vector(0,-1){1}}
\put(2,2){\vector(0,1){1}}

\put(5,0){figure 4.9}

\end{picture}\vspace{0.5cm}

d) This follows most of the time directly from duality. In the exceptional case $\kappa_{1}\alpha\gamma\lambda\gamma\neq 0$ an arrow $\phi\neq\kappa_{1}$ stopping at the end of $\kappa_{1}$ induces the cleaving diagram shown in the second picture of  figure 4.9.

\subsection{The proof of theorem 3}

 Throughout this section $P$ denotes a minimal representation-infinite ray category. Then there is no non-zero projective-injective $P$-module by a classical argument.

 Let $C$ be a non-deep contour in $P$. Then $C$ is a mild contour and we can apply the structure theorem and look at the three cases. If $C$ is a diamond then the projective $P(x,-)$ is also injective by proposition \ref{MID} which is impossible.

Next, let $C$ be a dumb-bell. Since $C$ is representation-finite and $P$ minimal representation-infinte there is  an arrow $\tau$ connected with $x$ or $y$ and different from $\lambda,\mu$ and $\rho$. Up to duality we get from lemma \ref{LDB} that $\tau$ starts in $y$ and we can apply proposition \ref{MIDB}.
If $\tau\mu=0$ then again $P(x,-)$ is projective-injective. If $\tau\mu \neq 0$ then we can split $z$ by proposition \ref{MIDB} into an emitter and a receiver. The obtained quiver has two connected components and the  indecomposables live only in one of these. In particular, there is no faithful indecomposable as it should be for a minimal representation-infinite category.

Finally, if $C$ is a penny-farthing with $P(x_{0},y) = 0 = P(y,x_{0})$ for all $y$ outside $C$ then $P(x_{0}, -)$ is again projective-injective. So up to duality we can assume $P(x_{0},b)\neq 0$ for some $b$ outside $C$. Then we use proposition \ref{MIPF} and the arguments from section 3 in \cite{Standard} to exclude this case.

\subsection{Another proof of theorem 3}

Let $P$ be a minimal representation-infinite ray category containing a non-deep contour $C=\{v,w \}$ from $x$ to $y$ such that $\vec{v}$ is transit. By theorem 2 in \cite{gaps} the universal cover $\tilde{P}$ of $P$ is interval finite and the fundamental group is free. By the finiteness criterion $\tilde{P}$ is not locally representation-finite. 

Assume that all finite convex subcategories of $\tilde{P}$ are representation-finite. Then $\tilde{P}$ contains arbitrarily long zigzags as convex subcategories that give rise to a crown in $P$. These zigzags involve no commutativity relations and so the crown exists already in the proper quotient  $P/\pi(y)\vec{v}$ contradicting the fact that $P$ is minimal representation-infinite.

\vspace{0.3cm}
\setlength{\unitlength}{0.8cm}
\begin{picture}(20,4)

\put(5,2){\vector(1,0){1}}
\put(6,2){\vector(1,0){1}}
\put(5,3){\vector(1,0){1}}
\put(6,3){\vector(1,0){1}}
\put(5,3){\vector(0,-1){1}}
\put(6,3){\vector(0,-1){1}}
\put(7,3){\vector(0,-1){1}}

\put(5,0){figure 4.10}
\end{picture}\vspace{0.5cm}

So there is a finite convex subcategory $B$ of $\tilde{P}$ whose frame belongs to the BHV-list. Since $P$ is minimal representation-infinite, $B$ cannot be annihilated by all the liftings of the path $p_{1}\ldots p_{r}v_{1}\ldots v_{n}$ introduced shortly before lemma \ref{wegeindeutig}. Thus the frame of $B$ has to contain the commutative diagram given in figure 1.2 as a subdiagram and so in particular two commutative $p$-gons for some $p$. A glance at the BHV-list reveals that only the frames 85,76,87 and 88 are possible. In particular we have $p=2$ and the two commutative squares are situated as in figure 4.10. This implies directly that the contour is a dumb-bell or a penny-farthing with $n=1$. Using the propositions \ref{MIDB} and \ref{MIPF} one completes the proof. In fact, one needs only very special cases of these propositions and the proofs thereof are much easier when one works in $\tilde{P}$.\vspace{1cm}

\subsection{Concluding remarks}
Note that a similar argument proves the somewhat surprising fact that a minimal representation-infinite algebra is defined by zero-relations and at most three commutativity relations. 

Also one knows that the universal cover $\tilde{P}$ is interval-finite for all mild categories. Thus a contour $C$ downstairs in $P$ can be lifted to a contour $\tilde{C}$ upstairs in $\tilde{P}$ and one might try to classify the possible non-deep contours that way. But this is a bad idea, because the structure of the faithful indecomposables over simply connected algebras is terribly complicated as is well-known from the second unpublished part of my habilitation. The zoo of these animals can be regarded in \cite{minimale}.  However, looking at the very simple list of the large indecomposables given in \cite{Treue}  one sees immediately that only penny-farthings can occur as large non-deep contours.

\end{document}